%% file: paper_v24.tex
\documentclass[%
3p,%
preprint,
fleqn, 
times
]{elsarticle}
\RequirePackage{fix-cm}
\usepackage{graphicx}
\usepackage{amsmath, amssymb}
\usepackage{amsthm, thmtools}
\usepackage{xcolor}
\usepackage{booktabs}
\usepackage{hyperref}
  \hypersetup{breaklinks=true,colorlinks=true,pdfborder={0 0 0}}
\newcounter{ARcount}[section]

\declaretheorem[style=definition, sibling=ARcount]{Definition}

\declaretheorem[style=definition, sibling=ARcount]{Remark}
\declaretheorem[style=definition, qed=$\blacksquare$, sibling=ARcount]{Theorem}
\declaretheorem[style=definition, sibling=ARcount]{Lemma}

\definecolor{golden}{rgb}{0.97, 0.83, 0.0}

\usepackage{stmaryrd} 
\allowdisplaybreaks
\input{Commands}

\usepackage{cancel}

\begin{document}

\journal{Computers \& Mathematics with Applications}
\title{Discontinuous Galerkin method for coupling hydrostatic free surface flows to saturated subsurface systems}

\author[FAU]{Andreas Rupp}
  \ead{rupp@math.fau.de}
\author[AWI,FAU]{Vadym Aizinger\corref{cor}}
  \ead{vadym.aizinger@awi.de}
\author[FAU]{Balthasar Reuter}
  \ead{reuter@math.fau.de}
\author[FAU]{Peter Knabner}
  \ead{knabner@math.fau.de}
\address[FAU]{Friedrich--Alexander University of Erlangen--N\"urnberg, Department of Mathematics, 
Cauerstra{\ss}e~11, 91058~Erlangen, Germany}
\cortext[cor]{Corresponding author}
\address[AWI]{Alfred Wegener Institute, Helmholtz Centre for Polar and Marine Research, Am Handelshafen 12, 27570 Bremerhaven, Germany}

\date{Received: date / Accepted: date}

%
\begin{abstract}
We formulate a~coupled surface/subsurface flow model that relies on hydrostatic equations with free surface in the free flow domain and on the Darcy model in the subsurface part. The model is discretized using the local discontinuous Galerkin method, and a~statement of discrete energy stability is proved for the fully non-linear coupled system.
\end{abstract}
\maketitle
\section{Introduction}\label{sec:introduction}
The interaction between free flow and subsurface systems (the latter either saturated or unsaturated) is important for a~variety of environmental applications, e.g. infiltration of overland flow into the soil during rainfall, contaminant propagation into the subsurface, sedimentation processes, interaction of seas, lakes, rivers, or wetlands with groundwater aquifers. Mathematical models for such coupled surface/subsurface flows generally express the conservation of mass and momentum in the coupled system. Coupled models usually pose substantial challenges on various levels: Mathematical -- due to differences in PDE system types in different subdomains giving rise to well-posedness and stability issues, numerical -- due to a~pronouncedly multi-scale character of the flow, and computational -- arising from the growing algorithmic complexity and increased performance and parallel scalability demands.

Depending on the target application and the level of modeling complexity, different model combinations in the surface and subsurface subdomains have been considered in the literature; the aspects covered include:
\begin{itemize}
\setlength\itemsep{0pt}
\item modeling approaches, in particular various choices of conditions at the model interface,
\item numerical methodology focusing on sub-problem discretizations and on solution algorithms that become critically important in the case of time-dependent flows,
\item theoretical issues mainly investigating the well-posedness, stability, and accuracy of coupled formulations.
\end{itemize}
In the context of geophysical flows, one can distinguish between two main types of fluid in the free flow subdomain, water and air, although both can certainly transport various additional substances. The subsurface systems usually contain either one (water), two (water and air), or three (water, air, and, e.g. oil) distinct phases, and each of those can furthermore transport additional species. The coupled models investigating the air flows (or transport of other gases dissolved in the air, e.g. CO$_2$) usually consider Stokes model in the surface subdomain and the one- or two-phase Darcy or Richards equations in the subsurface part~\cite{Rybak_etal_15,Mosthaf_Baber_etal_11}.

The modeling efforts for flow and transport processes involving water -- such as the present study -- cover a~much greater range of models in the free surface flow subdomain. A~number of recent studies (see \cite{Spanoudaki2009,Maxwell2014} for an~intercomparison) consider coupling free surface flows represented either by the 1D/2D shallow water equations~\cite{Dawson_08} or even simpler models (e.g. the kinematic wave equation~\cite{Sochala_Ern_Piperno_09}, a~diffusion wave approximation of the Saint--Venant equation~\cite{Sulis_etal_10}) with saturated subsurface flow described by the Richards or Darcy equations.

The theoretical aspects of coupled surface/subsurface flow modeling such as the well-posedness and the stability have also attracted some attention in the last decade. The most relevant studies in the context of the present work consider a~3D Navier-Stokes/Darcy-coupling based on a~discontinuous Galerkin (DG) method or on various combinations~\cite{Chidyagwai2009,Girault2013,Cesmelioglu2013,Badea2010} of the DG and finite element methods (see overview in~\cite{Discacciati2009}).

The hydrostatic primitive equations (sometimes also called the 3D shallow water equations) employed in our work is the most commonly used model for simulating circulation in geophysical domains with free surface such as oceans, lakes, estuaries, etc. The main assumption underlying this model (and setting it apart from the incompressible Navier-Stokes equations it is derived from) is the ratio between the horizontal and the vertical dimensions of at least 20:1 \cite[Sec.~2.3]{Vreugdenhil} with similar ratios for the horizontal to the vertical velocities and accelerations. This clear separation of the horizontal from the vertical scales is a~critical aspect of the hydrostatic modeling and is reflected in the direction of the gravity force, turbulence parametrizations, computational meshes made up of thin long elements with strictly vertical lateral faces, and many other details. This system also serves as the starting point for the derivation of the well known 2D shallow water equations. 

Although the hydrostatic primitive equations is a~widely used model, the aspects of well-posedness and stability of this PDE system as well as similar investigations of its discretizations are not very common and certainly appear to be neglected compared to more general models such as incompressible Navier--Stokes equations or less general ones such as 2D shallow water equations. The exceptions include works by Lions et al.~\cite{Lions1992a,Lions1992b}, Azerad~\cite{Azerad2001}, and the existence proofs for global strong solutions~\cite{Kobelkov2006,Titi2012}. Regarding the finite element analysis, one can note several works of Guill{\'e}n-Gonz{\'a}lez and co-workers treating this discretized system and its analysis as the limiting case of the Stokes system~\cite{GG2005,GG2015a,GG2015b} and our previous study of the DG method~\cite{AizingerPaper}. However, very few authors consider the problem in its full complexity and include the non-linear advection, free surface, or attempt to handle the difficulties arising from the hydrostatic approximation of the vertical velocity component. All aforementioned works except for~\cite{AizingerPaper} make the rigid lid assumption, \cite{GG2015b} introduces a~viscosity term into the continuity equation; other common simplifications include omitting the non-linear advection~\cite{GG2015a,GG2015b} and factoring out the vertical velocity~\cite{GG2005,Titi2012}.

The area of numerical modeling for subsurface applications in all its facets enjoys far more attention; this concerns the development and testing of new discretization techniques as well as their analysis. We refer the interested reader to a~recent article~\cite{DiPietro2014} for an overview.

The present study formulates a~coupled model consisting of the free surface flows represented by the three-dimensional hydrostatic equations and a~subsurface flow system modeled by Darcy's law. A~coupling condition is introduced based on a~special form of dynamic pressure, this coupling is then motivated using the weak formulation of the coupled system. The model equations are discretized using the local discontinuous Galerkin (LDG) method introduced in~\cite{DawsonAizinger2005} and further developed in~\cite{AizingerPDPN2013} for the hydrostatic free surface system and in~\cite{AizingerRSK2016,RuppKnabner2017,RuppKnabnerDawson2018} for Darcy's law. Finally, a~statement of semi-discrete energy stability is proved for the full non-linear formulation that also accounts for the dynamic free surface in the free flow domain.

The rest of the current paper is structured as follows. The next section introduces the mathematical models for the free surface and subsurface flow systems and proposes the interface conditions. In Sec.~\ref{sec:weak}, the weak problem formulation for the coupled problem is provided, and our choice of interface conditions is motivated by proving a~statement of weak steady-state stability for homogeneous boundary conditions. In Sec.~\ref{sec:discrete}, both problems are discretized using the LDG method, and a~statement of discrete stability is proved in Sec.~\ref{sec:analysis}. A~convergence study using the proposed formulation is given in Sec.~\ref{sec:numerical}, and a~conclusions section wraps up this work.
\section{Mathematical model}\label{sec:model}
\subsection{Computational domain}
A very important feature of 3D geophysical flow models is their natural an\-isotropy due to the gravity force acting in the vertical direction. This fact is usually reflected in the mathematical and numerical formulations as well as in the construction of computational domains and grids. The top boundary of most 3D surface flow domains is a dynamically changing surface whose movements correspond to time variations in the free surface elevation, although some models make a 'rigid lid' assumption to avoid increased computational costs connected with dynamically changing meshes.

\begin{figure}[ht!]
 \centering
 \includegraphics[width=\textwidth]{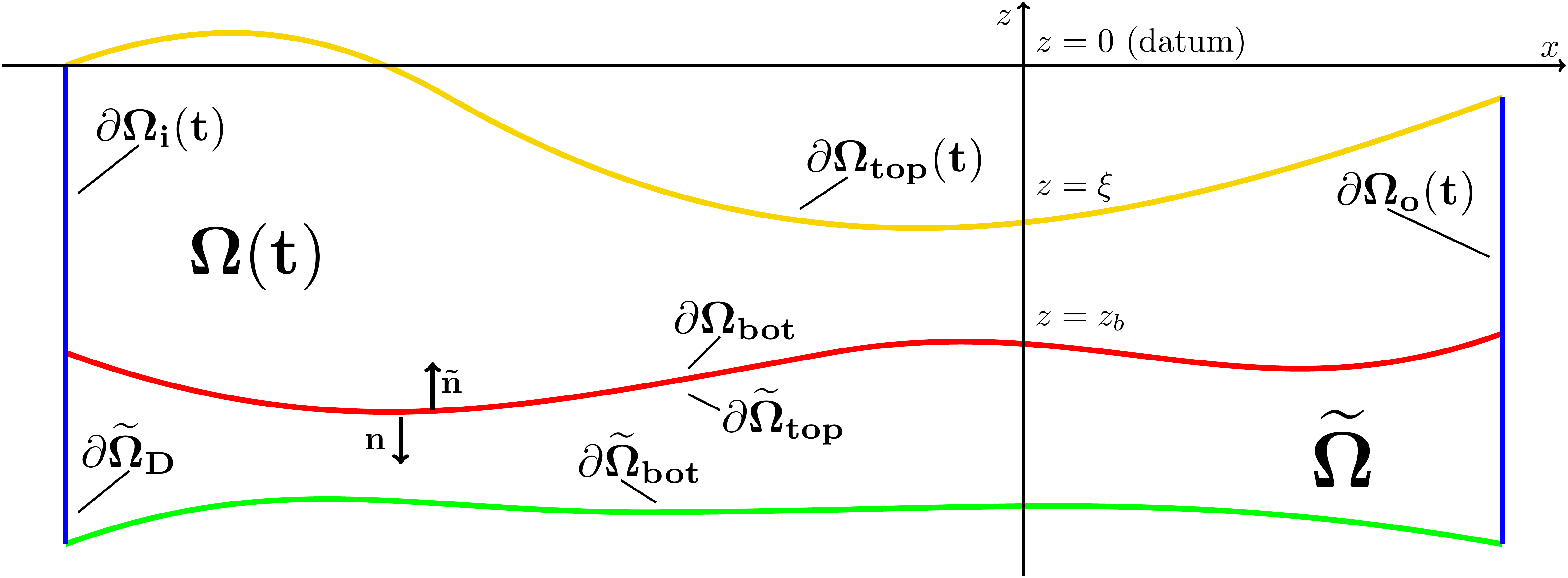}
 \caption{Vertical cross section of computational domain $\Omega(t)$ for hydrostatic equations on top of (fixed) computational domain $\tO$ for Darcy flow.} \label{figure0}
\end{figure}

Let $\Omega(t) \subset {\bR}^3$ (see Figure~\ref{figure0}) be our time-dependent domain for the hydrostatic free flow equations. We define $\Pi$ as the standard orthogonal projection operator from $\bR^3$ to $\bR^2$ ($\Pi(x, y, z) = (x, y)$, $\forall (x, y, z) \in \bR^3$), and $\Ox \coloneqq \Pi(\Omega(t))$. We require our top and bottom boundaries to be single-valued functions defined on $\Ox$ at any time (this excludes, e.g., wave breaking situations). 
The \textcolor{golden}{golden} top boundary of the domain $\p \Omega_{top}(t)$ is assumed to be the only moving boundary. The \textcolor{red}{red} bottom $\p \Omega_{bot}$ and \textcolor{blue}{blue} lateral $\p \Omega_{lat}(t)$ boundaries are considered to be fixed (though the height of the lateral boundaries can vary with time according to the movements of the free surface). We also require the lateral boundaries to be strictly vertical.
$\p \Omega_{bot}$ separates the time-dependent domain $\Omega(t)$ of the hydrostatic equations from the fixed domain $\tO$ of Darcy flow. Here, the \textcolor{red}{red} $\p \Omega_{bot} = \p \tO_{top}$ -- i.e. the bottom boundary of the free surface flow domain is the top boundary of the Darcy domain. The \textcolor{green}{green} boundary is the bottom boundary of the Darcy domain $\p \tO_{bot}$. 
In the following, all 2D counterparts of 3D vectors, operators, etc. consisting of the first two components of the former will be denoted by the subscript '2' without separate definitions (e.g., $\nablax \coloneqq (\p_x, \p_y)^T$). In a~similar manner, all functions defined on domain $\Ox$ will be trivially evaluated on $\Omega(t)$ via a~composition with $\Pi$, i.e. $\xi(x,y,z) \coloneqq \xi(\Pi(x,y,z))$. Furthermore, all unknowns, sets, etc. associated with the Darcy domain will be marked by tilde~$\tilde \cdot$.

\subsection{Primitive hydrostatic equations}\label{sec:hydrostatic}
The primitive hydrostatic equations with constant (unit) density describe the following properties of the free surface flow system \cite{Vreugdenhil}:
\begin{itemize}
 \item 2D conservation of volume (mass) also known as the primitive continuity equation (PCE)
 \begin{equation}\label{pce}
  \p_t \xi \ + \ \nablax \cdot \int_{z_b}^\xi \bux \, dz \ = \ F_H, 
 \end{equation}
 where $\xi, z_b$ are the values of the $z$ coordinate with respect to some datum at the free surface and the surface/subsurface flow interface, respectively, $\bu = (u,v,w)^T$ is the velocity vector, and $F_H$ is the source term that accounts for the normal flux from/to the subsurface domain. 
 \item 3D conservation of momentum (in conservative form)
 \begin{equation}\label{momentum_cons}
  \p_t \bux \ + \ \nabla \cdot \left(\bux \otimes \bu \ - \ \D \nabla \bux \right) \ + \ g \nablax \xi \ - \ \left(\begin{array}{cc} 0 & -f_c\\ f_c & 0 \end{array}\right) \bux \ = \ \bF_U,
 \end{equation}
 where the wind stress, the atmospheric pressure gradient, and the tidal potential are combined into a body force term $\bF_U$, $f_c$ is the Coriolis coefficient, $g$ is acceleration due to gravity. To prevent our analysis from being obscured by nonessential details, we simplify the momentum equations by omitting the Coriolis term and by rescaling the system so that $g=1$. The omission of the Coriolis term does not, in fact, affect the final result at all since the Coriolis force is energy neutral both, in the continuous and in the discrete sense and thus cancels out in the energy norm (see \cite{AizingerDiss}). In~\eqref{momentum_cons}, $\D = \D(\bu)$ denotes the tensor of eddy viscosity coefficients that can depend on the flow velocity (see, e.g., \cite{Davies1986}) defined as
 \begin{equation*}
  \D = \left( \ba{cc} D_u & 0\\ 0 & D_v \ea \right), \qquad \D \nabla \bux : = \left( \ba{c} (D_u \nabla u)^T\\ (D_v \nabla v)^T \ea \right) \in \IR^{2\times 3},
 \end{equation*}
 where $D_u$, $D_v$ and their inverses are $3 \times 3$ uniformly s.p.d. (symmetric positive definite) matrices.
 In Eq.~\eqref{momentum_cons}, $\otimes$ denotes the \emph{tensor product}, and $\nabla \cdot$ is the \emph{matrix divergence} defined as
 \begin{equation*}
  \left(\nabla \cdot A\right)_i \coloneqq \left(\sum_{j = 1}^n \p_{x_j} (A)_{i,j}\right)_i \qquad \text{for } i = 1, \ldots, m, \; A \in \IR^{m \times n}.
 \end{equation*}
 $D_u(\cdot, \cdot)$, $D_u^{-1}(\cdot, \cdot)$ being uniformly s.p.d. is equivalent to the existence of a constant $C_{D} \ge 1$ (independent of $\bu$) such that for all $\gvec x \in \IR^3$
 \begin{equation*}\label{REM:spdTens}
  C_{D}^{-1} \| \gvec x \|^2_2 \le \gvec x^T D_u(\bu)\, \gvec x \le C_{D} \| \gvec x \|^2_2, \qquad C_{D}^{-1} \| \gvec x \|^2_2 \le \gvec x^T D_u^{-1}(\bu)\, \gvec x \le C_{D}\, \| \gvec x \|^2_2.
 \end{equation*}
 This implies $C_{D} \ge \max\{\| D_u \|_{L^\infty(\Omega(t))}, \| D_u^{-1} \|_{L^\infty(\Omega(t))} \}$. For simplicity, we also assume that $D_v, D_v^{-1}$ satisfies the above equalities with the same constant $C_{D}$.
In the LDG framework employed in the present work, an~auxiliary variable $\q$ is introduced, and the second-order momentum equations \rf{momentum_cons} are re-written in mixed form
\begin{align}
 \p_t \bux \ + \ \nabla \cdot \left(\bux \otimes \bu \ + \q \right) \ + \nablax \xi & = \ \bF_U,\label{eq:mixed_momentum_1}\\
 \D^{-1}(\bu)\,\q + \nabla \bux & = 0,\label{eq:mixed_momentum_2}
\end{align}
where $\nabla\bux$ denotes the \emph{Jacobian} of $\bux$.
Note that Eqs.~\eqref{eq:mixed_momentum_1},\eqref{eq:mixed_momentum_2} actually represent a~system of $2+2\times3$ equations. 
\item 3D conservation of volume (mass) also known as the continuity equation
 \begin{equation}\label{cont_eq}
  \nabla \cdot \bu \ = \ 0.
 \end{equation}
Note that, differently from the incompressible Navier-Stokes system, \rf{cont_eq} is not a constraint used to determine pressure but rather an~equation for $w$.
\end{itemize}
The following boundary conditions (see \cite{Vreugdenhil} for details) are specified for the free surface flow system (except for the interface boundary given in Sec.~\ref{sec:interface}):
%
\begin{itemize}
 \item Denoting by $\bn = (n_x, n_y, n_z)^T$ an~exterior unit normal to the boundary of $\Omega(t)$ we distinguish between lateral inflow $\p \Omega_{i}(t) : = \{\p \Omega_{lat}(t) : \bu \cdot \bn \le 0\}$ and lateral outflow $\p \Omega_{o}(t) \coloneqq \p \Omega_{lat}(t) \setminus \p \Omega_{i}(t)$ boundaries. 
 \begin{equation}\label{hydrostatic_lateral_bc}
  {\bux}\big|_{\p \Omega_{lat}}\ = \ \buxh, \qquad \xi\big|_{\Pi(\p \Omega_{i}(t))} \ = \ \hat{\xi}.
 \end{equation}
 Even though the velocity is specified on the whole lateral boundary of $\Ot$, the advection terms only use the normal flux boundary condition at the outflow boundary $\buxh \cdot \bnx$ (see \eqref{discrete_general_2}). This somewhat unusual placement of flux (at $\Omega_o$) and water elevation (at $\Omega_i$) allows to compactify our discrete stability analysis and can be reversed to a~more standard configuration -- at the cost of some additional technicalities.
 \item The free surface boundary conditions have the form
 \begin{equation}\label{hydrostatic_surface_bc}
  \nabla u(\xi) \cdot \bn \ = \ \nabla v(\xi) \cdot \bn \ = \ 0.
 \end{equation}
 \item Additionally, initial data for $\bu$ and $\xi$ is given. Note that the initial and boundary conditions must be {\em compatible}.
\end{itemize}
%

%
%
Thus the free flow system that we consider in this problem consists of Eqs.~\eqref{pce}, \eqref{eq:mixed_momentum_1}, \eqref{eq:mixed_momentum_2}, \eqref{cont_eq} complemented by the corresponding initial and boundary conditions. Also note that the introduced simplifications neither affect the non-linearity of the system nor lower any analysis hurdles.
\subsection{System of equations for 3D Darcy flow}\label{sec:darcy}
Single phase flow through a porous medium $\tO$ is usually modeled by Darcy's law linking the hydraulic head $\tH$ and the seepage velocity $\tU = (\tilde u, \tilde v, \tilde w)$. In mixed formulation, the equations for constant (unit) density have the form:
\begin{subequations}
\begin{align}
 \p_t \tH + \nabla \cdot \tU & = \tilde f, \label{EQ:Darcy:cont1}\\
 \tD^{-1}(\tH) \tU + \nabla \tH & = 0 \label{EQ:Darcy:cont2}
\end{align}
for~a given source term $\tilde f$ and with $\tD$ and its inverse uniformly s.p.d. tensors (similarly to Sec.~\ref{sec:hydrostatic} and with the same constant $C_{D}$). The boundary conditions for the flux and head are given by
\begin{equation}\label{darcy_bc}
 (\tU \cdot \tilde \bn)\big|_{\p \tO_N}  = \hat u_{\tilde n}, \qquad \tH{\big|_{\p \tO_D}}  = \hat h.
\end{equation}
\end{subequations}
Here, $\tilde \bn$ denotes the outward unit normal with respect to $\tO$. In addition to this, initial data $\tH_0$ is given. The bottom and lateral boundaries of $\tO$ are either Dirichlet or Neumann boundaries, while a coupling boundary condition is imposed at the top boundary. Eqs.~\eqref{EQ:Darcy:cont1},\eqref{EQ:Darcy:cont2} have been simplified via division by the specific storativity $\tilde \Phi(t, \vec x) \ge \tilde \Phi_0 > 0$.
%
\subsection{Interface conditions}\label{sec:interface}
Specifying the interface conditions between the sub-models is not a simple task in the context of the present study; the main difficulty is finding a~set of transition conditions that guarantee a~physically founded and mathematically well-posed system of equations for the coupled model. In our case, this task is more challenging for the free flow model due to its greater complexity (i.e., the presence of non-linear advection terms). Thus, even a~standard variational formulation of the incompressible Navier-Stokes/Darcy system includes an undetermined-sign term (see \cite{Girault2009}) on the transition boundary. By resorting to a~linear Stokes model some authors avoid this problem (see discussions of the modeling and coupling issues in \cite{Discacciati2009,Girault2013}). Another avenue to handle this problem involves modifying the momentum equations by adding the so called Temam stabilization \cite{Temam1968} term that is equal to zero in the strong sense but can be exploited in a~way that provides some additional control over the kinetic energy in the weak formulation. 

In this study, no Temam stabilization is used, and the full non-linear advection is retained. We impose the following transition conditions at the boundary $\p \Ot_{bot}$ = $\p \tO_{top}$ between the free surface and subsurface flow subdomains:
\begin{itemize}
 \item the continuity of the normal flux (volume/mass conservation) 
 \begin{equation}\label{transition_bc_1}
  (\bu \cdot \bn)\big|_{\p \Ot_{bot}} = - (\tU \cdot \tilde \bn)\big|_{\p \tO_{top}};
 \end{equation}
 \item the continuity of pressure (head), where we use a~special form of dynamic pressure in the free flow subdomain (cf.~\cite{Girault2009,Fetzer2016}) and ignore viscous terms (also see the discussion in Sec.~\ref{sec:weak-darcy})
 \begin{equation}\label{transition_bc_2}
  \tH\big|_{\p \tO_{top}} = \left( \xi + \frac{\bux \cdot \bux}{2 g}\right)\bigg|_{\p \Ot_{bot}} \quad \mbox{($g$=1 was assumed in Sec.~\ref{sec:hydrostatic} and is included here for consistency)},
 \end{equation}
where we recall the hydraulic head definition: $\tH = z + p/(g \rho_w)$, with $z$ denoting the vertical coordinate of the point with respect to the datum, $\rho_w$ the water density, and $p$ the fluid pressure;
 \item the friction law on horizontal velocity components modeled on the standard friction laws for turbulent shallow-water flows (see, e.g., \cite{Vreugdenhil}) and rather similar to the Beavers-Joseph-Saffman~\cite{Saffman1971} condition very common in coupled surface/subsurface flow applications
 \begin{equation}\label{transition_bc_3}
  \D_u \nabla u(z_b) \cdot \bn \ = \ -C_f(\bu) u(z_b), \qquad \D_v \nabla v(z_b) \cdot \bn \ = \ -C_f(\bu) v(z_b),
 \end{equation}
 where the minus sign is due to $\bn$ being an {\em exterior} unit normal to $\p \Ot_{bot}$, and $C_f(\bu)>0$ is the bottom friction coefficient that in shallow-water applications is usually represented by either $C_f(\bu)=const$ for a~linear or by $C_f(\bu) = C'_f \,|\bu(z_b)|$ with $C'_f = const$ for a~quadratic friction law.
\end{itemize}
The interface conditions specified above are modeled closely on those used in Navier-Stokes/Darcy coupled models (these are the closest analog to our setting found in the literature, see, e.g.~\cite{Fetzer2016}) with certain modifications motivated by the important differences between the incompressible Navier-Stokes and the hydrostatic model used here. The main difference reflected both in the dynamic pressure term \eqref{transition_bc_2} and in the friction formula \eqref{transition_bc_3} is the fact that a~hydrostatic system does not conserve the vertical momentum; instead, the vertical velocity is computed by the continuity equation~\eqref{cont_eq} that expresses the 3D conservation of mass/volume. This circumstance makes a~physically consistent formulation of coupling conditions for the momentum equations particularly challenging. 
\section{Weak formulation of the coupled system}\label{sec:weak}
\subsection{Weak formulation of the hydrostatic equations}\label{sec:weak-swe}
To simplify notation we use from now on $\| u \|_\Omega$ for the $L^2$ norm of $u$ and $( \ . \ , \ . \ )_\Omega$, $< \ . \ , \ . \ >_\gamma$ for the $L^2$ inner products on domains $\Omega \subset {\bR}^d$ and surfaces $\gamma \subset {\bR}^d$, respectively. Used in conjunction with vectors or tensors, these products are to be understood as sums of componentwise $L^2$ inner products.

Next, we obtain a~weak formulation of the hydrostatic system by multiplying Eqs.~\eqref{pce}, \eqref{eq:mixed_momentum_1}, \eqref{eq:mixed_momentum_2}, \eqref{cont_eq} with some smooth test functions and integrating by parts. For the PCE, we get:
\begin{equation*}
\left(\p_t \xi, \delta \right)_\Ox 
+ \lan \int_{z_b}^\xi \bux \, dz \cdot \bnx, \delta \ran_{\p \Ox} 
\hspace{-2mm} - \left( \int_{z_b}^\xi \bux \, dz \cdot \nablax, \delta \right)_\Ox  
= \left( F_H, \delta \right)_\Ox.
\end{equation*}
Exploiting the fact that the lateral boundaries of $\Ot$ are strictly vertical and substituting \eqref{transition_bc_1} into $F_H$, we can rewrite the equation above in a special 2D/3D form
\begin{subequations}\label{eq:weak-hydrostatic}
\begin{equation}\label{eq:weak-h}
 \left(\p_t \xi, \delta \right)_\Ox + \lan \bux \cdot \bnx, \delta \ran_{\p \Omega_{lat}} - \left(\bux \cdot \nablax, \delta \right)_\Ot + \lan \tilde \bu \cdot \bn, \delta \ran_{\p \Omega_{bot}} = 0.
\end{equation}
Note that \eqref{eq:weak-h} is well defined for any $\xi, \delta \in H^1(\Ox)$ and a.e. $t \in [0, T]$.

A~weak form of the momentum equations given by 
%
\begin{align}
 & \left(\p_t \bux, \bphi \right)_\Ot + \lan (\bux \otimes \bu +  \q) \cdot \bn + \xi \bnx, \bphi \ran_{\p \Ot}
 - \left( (\bux \otimes \bu +  \q) \cdot \nabla + \xi \nablax, \bphi \right)_\Ot  = \left( \bF_U, \bphi \right)_\Ot,\label{eq:weak-u} \\
 & \left(\D^{-1}(\bu)\;\q, \bpsi \right)_\Ot + \lan \bux \otimes \bn, \bpsi \ran_{\p \Ot}\ - \ \left( \bux \otimes \nabla, \bpsi \right)_\Ot = 0.\label{eq:weak-q}
\end{align}
%
For the continuity equation, we get 
\begin{equation}\label{eq:weak-w}
 \lan \bu \cdot \bn, \sigma \ran_{\p \Ot}\ - \ \left(\bu \cdot \nabla, \sigma \right)_\Ot\ = \ 0.
\end{equation}
\end{subequations}
Eqs.~\eqref{eq:weak-u}--\eqref{eq:weak-w} are well defined $\forall \bux, \bphi \in H^1(\Ot)^2$, $\q, \bpsi \in H^1(\Ot)^{2 \times 3}$, $w, \sigma \in H^1(\Ot)$, and for a.e. $t \in [0, T]$. 
\subsection{Weak formulation of Darcy equation}\label{sec:weak-darcy}
For smooth test functions $\testH$ and $\testU$, we multiply \eqref{EQ:Darcy:cont1} and \eqref{EQ:Darcy:cont2} by these test functions and integrate by parts. 
\begin{subequations}\label{eq:weak-darcy}
\begin{align}
 & \left(\p_t \tH, \testH \right)_{\tO} 
- \left(\tU \cdot \nabla , \testH \right)_{\tO} 
+ \left\langle \tU \cdot \tilde \bn, \testH\right\rangle_{\p \tO} 
= \left(\tilde f, \testH \right)_{\tO},\label{eq:weak-tilde-h}\\
 & \left(\D^{-1}(\tH)\, \tU, \testU \right)_{\tO} 
- \left(\tH\, \nabla , \testU \right)_{\tO} 
+ \left\langle \tH \, \tilde \bn, \testU \right\rangle_{\p \tO} = 0\label{eq:weak-tilde-u}.
\end{align}
\end{subequations}
The above terms are well defined for $\tH,\testH \in H^1(\tO)$ and $\tU,\testU \in H^1(\tO)^3$, and for a.e. $t \in [0,T]$.
\subsection{Weak energy estimate for the coupled system}\label{sec:weak-coupled}
In this section, we formulate a~statement of weak energy stability for the coupled system to illustrate the difficulties connected with finding a~workable set of transition conditions at the coupling interface. We consider a~stationary variant of problem \eqref{eq:weak-hydrostatic}, \eqref{eq:weak-darcy} and further simplify our task by using homogeneous boundary conditions for velocities and fluxes in both, free flow and subsurface subdomains. That is $\partial \Omega_{lat} = \partial \Omega_i$ and $\partial \tO \setminus \partial \tO_{top} = \p \Omega_N$.

\noindent
Denoting by $H^1_{0,\Gamma}(\Omega)^d$ for $\Gamma \subset \p \Omega$ the space $\{f \in H^1(\Omega)^d: f\big|_{\Gamma} = 0\}$, we select the test and trial spaces as follows:
\begin{equation*}
\xi, \delta \in H^1(\Omega_2), \quad \bux, \bphi \in H^1_{0,\partial \Omega_{lat}}(\Omega)^2, \quad \q, \bpsi \in H^1_{0, \partial \Omega_{top}} (\Omega)^{2\times3}, \qquad
 \tH, \tilde \delta \in H^1(\tO), \qquad \tU, \tilde \bphi \in H^1_{0, \partial \tO \setminus \partial \tO_{top}}(\tO)^3.
\end{equation*}
Setting $\delta = \xi,\, \bphi = \bux,\, \bpsi = \q$ in \eqref{eq:weak-h}, \eqref{eq:weak-u}, \eqref{eq:weak-q} and using the definitions of test spaces and boundary conditions, we obtain
\begin{align*}
 & -(\bux, \nabla_2 \xi)_{\Ot} + \lan \tU \cdot \bn, \xi \ran_{\p \Omega_{bot}} = 0\\
 & \lan \bux (\bu \cdot \bn) + C_f \bux + \xi \bnx, \bux \ran_{\p \Omega_{bot}} + \lan \bux (\bu \cdot \bn) +  \xi \bnx, \bux \ran_{\p \Omega_{top}}
 - \left( \bux(\bu \cdot \nabla) +  \q \cdot \nabla + \xi \nablax, \bux \right)_\Ot = \left( \bF_U, \bux \right)_\Ot, \\
 & \left(\D^{-1}(\bu)\;\q, \q \right)_\Ot + \lan \bux, \q \cdot \bn \ran_{\p \Omega_{bot}} - \left( \bux , \nabla \cdot \q \right)_\Ot = 0.
\end{align*}
Since $\p_t \xi =0$, we have $(\bu \cdot \bn)\big|_{\p \Omega_{top}}=0$. Also note that the integration by parts and the continuity equation~\eqref{cont_eq} give us
\vspace{-2mm}
\begin{equation*}
 \left( \bux(\bu \cdot \nabla), \bux \right)_\Ot = \f 12 \left( \bu, \nabla |\bux|^2 \right)_\Ot = \f 12 \lan \bux (\tilde \bu \cdot \bn), \bux \ran_{\p \Omega_{bot}}.
\end{equation*}
%

Adding all equations together and using some simplifications that utilize the boundary conditions and an~integration by parts of element integral terms, we obtain the statement for energy in the free flow subdomain
\begin{align*}
 &\lan \tU \cdot \bn, \xi \ran_{\p \Omega_{bot}} + \lan C_f \bux, \bux \ran_{\p \Omega_{bot}} + \left(\D^{-1}(\bu)\;\q, \q \right)_\Ot 
 + \f 12 \lan \bux (\tilde \bu \cdot \bn), \bux \ran_{\p \Omega_{bot}} = \left( \bF_U, \bux \right)_\Ot.
\end{align*}
Setting $\testH = \tH, \testU = \tU$ in \eqref{eq:weak-tilde-h}--\eqref{eq:weak-tilde-u}, integrating by parts, adding equations, and using the boundary conditions, we get
\begin{equation*}
 \left(\tilde \D^{-1}(\tH)\, \tU, \tU \right)_{\tO} + \left\langle \xi + \f {\bux \cdot \bux} {2}, \tU \cdot \tilde \bn \right\rangle_{\p \tO_{top}} = \left(\tilde f, \tH \right)_{\tO}.
\end{equation*}
Since $\bn = - \tilde \bn$ on the interface boundary, the mass flux terms cancel out; also our choice of transition condition on the pressure becomes obvious. The statement of energy stability for the coupled system then reads:
\begin{equation*}
 \left(\D^{-1}(\bu)\, \q, \q \right)_\Ot + \lan C_f \bux, \bux \ran_{\p \Omega_{bot}} + \left(\tilde \D^{-1}(\tH)\, \tU, \tU \right)_{\tO} = \left( \bF_U, \bux \right)_\Ot +\left(\tilde f, \tH \right)_{\tO}.
\end{equation*}
\section{Discrete Formulation}\label{sec:discrete}
\subsection{Basic definitions and mathematical analysis tools}\label{sec:tools}
In the following, $\Th$ denotes a~non-overlapping $d$-dimensional polytopic partition of $\Omega \in \{\Ot, \tO, \Ox\}$ (see \cite[Def.~1.12]{PietroErn}). All partitions are assumed to be \emph{geometrically conformal} (in the sense of \cite[Def. 1.55]{ErnGuerm}). All proofs and arguments also hold for geometrically non-conformal meshes, but the notation becomes more cumbersome. The \emph{test} and \emph{trial} spaces for our LDG method are defined as the $d$-dimensional ($d \ge 1$) \emph{broken polynomial spaces of order $k$}
\begin{equation*}
 \IP_k^d(\Th) \coloneqq \left\{ \vec v \in L^2(\Omega)^d \, : \, \vec v_{\mid \elem} \mbox{ is a polynomial of degree at most $k$, } \forall \elem \in \Th \right\}.
\end{equation*}
Let $\F = \F(\Th)$ be the set of faces; for a~scalar function $w$ and a~vector function $\vec v$, we define the average $\avg{ \cdot }$ and the jump $\jump{ \cdot}$ on $\p \elem_i \cap \p \elem_j$ for neighboring mesh elements $\elem_i, \elem_j \in \Th,\, \elem_i \neq \elem_j$ in the following way:
\begin{align*}
 &\avg{ w }  = \frac{1}{2} \left( w_{\mid \elem_i} + w_{\mid \elem_j} \right), 
 &&\avg{ \vec v} = \frac{1}{2} \left( \vec v_{\mid \elem_i} + \vec v_{\mid \elem_j} \right),\\
 &\jump{ w } = w_{\mid \elem_i}\bn_{\elem_i} + w_{\mid \elem_j}\bn_{\elem_j}, 
 && \jump{ \vec v } = \vec v_{\mid \elem_i} \cdot \bn_{\elem_i} + \vec v_{\mid \elem_j} \cdot \bn_{\elem_j},
 &\jumpt{ \vec v } = \vec v_{\mid \elem_i} \otimes \bn_{\elem_i} + \vec v_{\mid \elem_j} \otimes \bn_{\elem_j},
\end{align*}
where $\bn_\elem$ is the outward unit normal with respect to $\elem$. Note, that a jump in a scalar variable is a vector, whereas a jump of a vector is a scalar. In addition, $\jumpt{\cdot}$ is introduced for vectors to denote a~second order tensor resulting from using the scalar jump definition component-wise. 
In our analysis, we use some well-known properties of jumps:
\begin{subequations}
\begin{align}
 \jump{ab} & = \; \avg{a} \jump{b} + \jump{a} \avg{b},\label{jump1}\\
 \avg{ab} & = \; \avg{a} \avg{b} + \frac 14 \jump{a} \cdot \jump{b}.\label{jump2}
\end{align}
\end{subequations}
The standard mathematical analysis tools used in this work include Young's and Cauchy-Schwarz' inequalities as well as the following results (see \cite[Sec.~1.4.1--1.4.3]{PietroErn})
\begin{Definition}[Shape and contact regularity]\label{def:regularity}
 A family of meshes $\Th$ is called \emph{shape and contact regular} (for short \emph{regular}) if, for all $\Delta x> 0$, $\Th$ admits a geometrically conformal, \emph{matching simplicial submesh} $\bar \Th$ such that 
 \begin{enumerate}
  \item $\bar \Th$ is \emph{shape-regular} in the sense of \cite{CiaHB}, i.e. there exists $\lambda_1 > 0$, independent of $\Delta x$, such that for all $\bar \elem \in \bar \Th$
  \begin{equation*}
   \lambda_1 \Delta x_{\bar \elem} \le  \rho_{\bar \elem},
  \end{equation*}
  where $\rho_{\bar \elem}$ is the diameter of the largest ball that can be inscribed in $\bar \elem$.
  \item there exists a~constant $\lambda_2 > 0$ independent of $\Delta x$ such that for all $\elem \in \Th$ and for all $\bar \elem \in \bar \Th$ with $\bar \elem \subset \elem$
  \begin{equation*}
   \lambda_2 \Delta x_\elem \le \Delta x_{\bar \elem},
  \end{equation*}
  \item there exists a~constant $\lambda_3 > 0$ independent of $\Delta x$ such that for all $\git \in \F$
  \begin{equation*}
   \lambda_3 \Delta x \le \Delta x_{\git}.
  \end{equation*}
 \end{enumerate}
\end{Definition}
\begin{Lemma}[Discrete trace inequality]\label{lem:disc:trace}
 Let $(\Th)$ be a regular mesh sequence with parameters $\lambda_1, \lambda_2, \lambda_3$. Then for all $\Delta x > 0$, all $\vec p \in \IP^d_k(\Th)$, the following holds with $C_t$ only depending on $\lambda_1, \lambda_2, \lambda_3$, $d$, and $k$:
 \begin{equation*}
  \Delta x^{1/2} \sum_{\git \in \F} \| \vec p \|_{L^2(\git)} \; \le \; C_t \sum_{\elem \in \Th} \| \vec p \|_{L^2(\elem)} \;=\; C_t \| \vec p \|_{L^2(\Omega)}.
 \end{equation*}
For $\git$ shared by elements $\elem_i$ and $\elem_j$, $\| \vec p \|_{L^2(\git)}$ is assumed to contain both traces
 \begin{equation*}
  \| \vec p \|_{L^2(\git)} \;=\; \| \vec p_{\mid \elem_i} \|_{L^2(\git)} + \| \vec p_{\mid \elem_j} \|_{L^2(\git)}.
 \end{equation*}
\end{Lemma}
%
%
\subsection{Computational mesh and free surface representation}
Keeping in line with the specific anisotropy of $\Omega(t)$ we construct our 3D mesh by extending a 2D triangular mesh of $\Ox$ in the vertical direction resulting in a~3D mesh of $\Omega(t)$ that consists of one or more layers of prismatic elements. In order to better reproduce the bathymetry and the free surface elevation of the computational domain, top and bottom faces of prisms can be non-parallel to the $xy$-plane; however, the lateral faces are assumed to be strictly vertical. 

\begin{figure}[h!]
\center
 \includegraphics[width=0.92\textwidth]{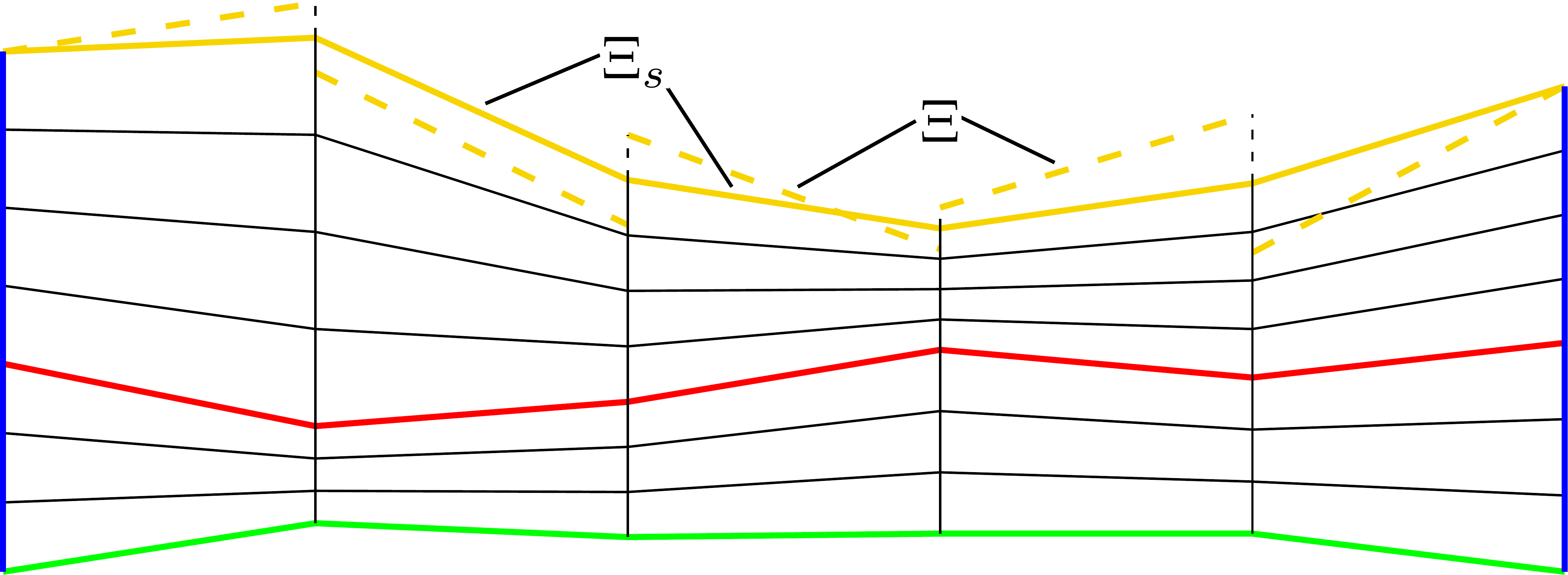}
 \caption{Vertical cross-section of the coupled mesh and the free surface geometry approximation (solid yellow line).}\label{Mesh_smoothing}
\end{figure}

For our analysis, we introduce the following sets of elements and faces:
\begin{itemize}
 \item $I_e$ - set of prismatic elements in $\Omega(t)$;
 \item $I_{e, 2D}$ - set of triangular elements in $\Ox$;
 \item $I_{e, \elemx}$ - set of prismatic elements corresponding to 2D element $\elemx$; 
 \item $I_{lat}$ - set of interior lateral faces in $\Omega(t)$; 
 \item $I_{horiz}$ - set of interior horizontal faces in $\Omega(t)$;
 \item $I_i, I_o$ - sets of exterior inflow and outflow lateral faces in $\Omega(t)$;
 \item $I_{top}$ - set of exterior faces on the top boundary of $\Omega(t)$;
 \item $I_{bot}$ - set of exterior faces on the bottom (transition) boundary of $\Omega(t)$;
 \item $\tIe$ - set of elements in $\tO$;
 \item $\tIi$ - set of interior faces in $\tO$;
 \item $\tID$ - set of faces on Dirichlet boundary of $\tO$;
 \item $\tIT$ - set of faces on top (transition) boundary of $\tO$;
 \item $\tIN$ - set of faces on Neumann boundary of $\tO$.
\end{itemize}
A key feature of our 3D LDG model is the fact that all primary variables -- including the free surface elevation -- are discretized using discontinuous polynomial spaces. As a result, computed values of the free surface elevation may have jumps across inter-element boundaries. If our finite element grids were to follow exactly the computed free surface elevation field this would cause the elements in the surface layer to have mismatching lateral faces (staircase boundary). We avoid this difficulty by employing a globally continuous (piecewise linear) free surface approximation that is obtained from the computed values of the free surface elevation with the help of a smoothing algorithm (see Fig.~\ref{Mesh_smoothing}) and denote by $\Xi_s$ the free surface elevation of the smoothed mesh. It must be noted here that solely the computational mesh is modified by the smoothing algorithm whereas the computed (discontinuous) approximations to all unknowns, including the free surface elevation, are left unchanged. This approach preserves the local conservation property of the LDG method and is essential for our algorithm's stability.
\subsection{Semi-discrete LDG formulation for the hydrostatic equations}
Our next step is to approximate $\left(\xi(t,\cdot), \bux(t,\cdot), w(t,\cdot), \q(t,\cdot)\right)$, a~solution to the weak problem, with a~function $\left(\Xi(t,\cdot), \bUx(t,\cdot), W(t,\cdot), \Q(t,\cdot)\right)\in {\cal H}_\Delta \times U_\Delta \times W_\Delta \times Z_\Delta$, where ${\cal H}_\Delta$, $U_\Delta$, $W_\Delta $, and $Z_\Delta$ denote finite-dimensional DG spaces.  
For this purpose, we use the weak formulation with one important modification: Since the DG approximation spaces do not guarantee continuity across the inter-element boundaries, all integrands in the integrals over interior faces have to be approximated by suitably chosen numerical fluxes that preserve consistency and stability of the method. Similar treatment may be needed at the exterior boundaries as well. Then a semi-discrete finite element solution is obtained by requiring that for a.e. $t \in [0, T]$, for all $\elem \in \Th$, and for all $(\delta,\bphi,\bpsi,\sigma) \in {\cal H}_\Delta \times U_\Delta \times W_\Delta \times Z_\Delta$, the following holds:
\begin{subequations}
 \begin{align}
& \hspace{-5mm} \left(\p_t \Xi, \delta \right)_\elemx + \sum_{\elem \in I_{e, \elemx}} \left\{ \lan  R_H, \delta \ran_{\p \elem_{lat}} - \left(\bUx \cdot \nablax, \delta \right)_\elem +  \lan \TU \cdot \bn , \delta \ran_{\p \elem \cap \p \Omega_{bot}} \right\}
= 0, \label{discrete_general_1}\\
& \hspace{-5mm}  \left(\p_t \bUx, \bphi \right)_\elem + \lan \mathbf{R}_U + \mathbf{S}_U, \bphi \ran_{\p \elem} - \left( (\bUx \otimes \bU +  \Q) \cdot \nabla + \Xi \nablax, \bphi \right)_\elem 
  + \f {n_z} 2 \lan \p_t (\Xi_s-\Xi)\, \bUx, \bphi \ran_{\p \elem \cap \p \Omega_{top}} = \left(\bF_U, \bphi \right)_\elem, \label{discrete_general_2}\\
& \hspace{-5mm}  \left(\D^{-1}(\bU)\;\Q, \bpsi \right)_\elem + \lan S_Q, \bpsi \ran_{\p \elem} - \left( \bUx \otimes \nabla, \bpsi \right)_\elem = 0 , \label{discrete_general_3}\\
& \hspace{-5mm}  \lan R_H, \sigma\ran_{\p \elem_{lat}} + \lan \bUd \cdot \bn, \sigma \ran_{\p \elem_{horiz} \setminus \p \Omega_{bot}} + \lan \TU \cdot \bn, \sigma \ran_{\p \elem \cap \p \Omega_{bot}} -  \left(\bU \cdot \nabla, \sigma \right)_\elem = 0,\label{discrete_general_4}
 \end{align}
\end{subequations}
where $\bUd$ denotes the value of $\bU$ taken from the element below the horizontal face, and $\p \elem_{lat}$ and $\p \elem_{horiz}$ are the lateral and the horizontal parts (faces) of $\p \elem$, respectively. $R_H$, and $\mathbf{R}_U$ are the normal advective fluxes for $\bUx \cdot \bnx$ and $\bUx (\bU \cdot \bn) + \Xi \bnx$ (also see Remark~\ref{remark-advective-fluxes}), respectively, whereas $\mathbf{S}_U$ and $S_Q$ denote the normal diffusive fluxes on element faces (Remark~\ref{remark-diffusive-fluxes}).
In our study, the following flux approximations are used:
\begin{subequations}
\begin{align}
\label{flux_s}
& \mathbf{S}_U|_{\gamma} 
\coloneqq \left\{\begin{array}{ll}
\avg{\Q} \cdot \bn, & \quad\gamma \in I_{lat} \cup I_{horiz},\\
\Q \cdot \bn, & \quad\gamma \in I_i \cup I_o,\\
0,& \quad\gamma \in I_{top},\\
C_f \bUx, & \quad\gamma \in I_{bot},
\end{array} \right.\quad
 S_Q|_{\gamma} 
\coloneqq \left\{\begin{array}{ll}
\avg{\bUx} \otimes \bn, & \quad\gamma \in I_{lat} \cup I_{horiz},\\
\buxh \otimes \bn, & \quad\gamma \in I_i \cup I_o,\\
\bUx \otimes \bn,& \quad\gamma \in I_{top} \cup I_{bot}.
\end{array} \right.\\
\label{flux_r}
&  R_H\big|_{\gamma} 
\coloneqq \left\{\begin{array}{ll}\avg{\bUx} \cdot \bnx,& \quad\gamma \in I_{lat},\\ 
\bUx \cdot \bnx, & \quad \gamma \in I_i,\\ 
\buxh \cdot \bnx, & \quad \gamma \in I_o, \end{array} \right. \quad
\mathbf{R}_U|_{\gamma} 
\coloneqq \left\{\begin{array}{ll} \avg{\bUx \otimes \bU} \cdot \bn + \avg{ \Xi } \bnx + \f {\lambda_U} 2 \jumpt{\bUx} \cdot \bnx,& \quad\gamma \in I_{lat},\\
\left(\avg{\bUx} \otimes \bUd\right) \cdot \bn + \Xi \bnx,& \quad\gamma \in I_{horiz},\\
\bUx (\bU \cdot \bn) + \hat{\xi} \bnx + \f {\lambda_U} 2 (\bUx - \buxh),& \quad\gamma \in I_i,\\
\bUx (\buh \cdot \bn) + \Xi \bnx,& \quad\gamma \in I_o,\\
\bUx (\bU \cdot \bn) + \Xi \bnx,& \quad\gamma \in I_{top},\\
\bUx (\TU \cdot \bn) + \Xi \bnx,& \quad\gamma \in I_{bot}
\end{array} \right.
\end{align}
\end{subequations}

The value of the penalty coefficient in the momentum flux on lateral faces is closely related to that of the standard Lax-Friedrichs solver (see a~discussion in the remainder of this section) and is given by
\begin{equation}
  \lambda_U \coloneqq \left\{\begin{array}{ll} \avg{ \left|\bUx \cdot \bnx\right| } + \sqrt{\avg{\left|\bUx \cdot \bnx\right|}^2 + 1}, & \mbox{ on interior lateral faces,}\\ \left|\bUx \cdot \bnx\right| + \sqrt{\left(\bUx \cdot \bnx\right)^2 + 1}, & \mbox{ on inflow lateral faces.} \end{array} \right.\label{lambda_U_def}
\end{equation}
Also note that the vertical component of the normal to lateral faces is zero, thus $\bU \cdot \bn = \bUx \cdot \bnx$, etc. there.

\begin{Remark}[Advective fluxes] \label{remark-advective-fluxes}
Normal fluxes $R_H$ and $\mathbf{R}_U$ for the non-linear advection operator on lateral faces for the PCE~\rf{discrete_general_1} and the momentum Eqs. \rf{discrete_general_2} must be computed by solving a~Riemann problem in a~coupled way (see \cite{AizingerDiss,DawsonAizinger2005,Toro2001} for a~discussion of this issue). These fluxes are much more important for the stability of the discrete scheme than those on horizontal faces. This phenomenon is a consequence of the specific anisotropy of our problem and of the computational mesh tailored for this anisotropy: the free surface elevation has jumps across lateral faces but not across horizontal ones. On horizontal faces, the free surface elevation $\Xi$ is continuous, thus the Riemann problem simplifies to that for momentum equations only. 
\end{Remark}

The largest (in absolute value) eigenvalue of the normal advective flux (see \cite{AizingerPaper}) given by
\begin{equation*}
 |\lambda_{\max} (\bUx)| = |\bUx \cdot \bnx| + \sqrt{ (\bUx \cdot \bnx)^2 + 1}
\end{equation*}
is used in the standard Lax-Friedrichs flux as the penalty coefficient. 
In this work, we slightly modify this Riemann solver to reduce the technicalities involved in the the stability analysis; however, the standard Lax-Friedrichs formulation works as well. The modifications amount to just omitting the penalty term in the PCE \eqref{discrete_general_1} and retaining it in the momentum \eqref{discrete_general_2} equation (cf.~\eqref{flux_r}).

For our choice of penalty coefficient given in Eq.~\eqref{lambda_U_def}, one can prove the following
\begin{Lemma}[Properties of $\lambda_U$\label{properties_lambda}]
 The following inequality holds for $\lambda_U$:
 \begin{equation*}
  \lambda_U(\bUx) \; \ge \; \f {\sqrt{2} + 1}{\sqrt{2}} \left| \bUx \cdot \bnx \right| + \f 1{\sqrt{2}}.
 \end{equation*}
\end{Lemma}
\noindent
{\em Proof}: This property follows directly from a simple arithmetic inequality
\begin{equation*}
 a + b \; \le \; \sqrt{2a^2 + 2b^2}, \qquad \forall a,b \ge 0.
\end{equation*}
\begin{Remark}[Diffusive fluxes]\label{remark-diffusive-fluxes}
Choosing diffusive fluxes $\mathbf{S}_U$ and $S_Q$ in Eqs.~\eqref{discrete_general_2} and \eqref{discrete_general_3} is simpler than solving the corresponding problem for the advective fluxes. In our analysis and implementation, those were set equal to central approximations on interior faces and to corresponding boundary conditions on the exterior ones (see~\eqref{flux_s}.
\end{Remark}

\begin{Remark}[Mesh penalty]\label{remark-mesh-penalty}
In addition to the ``usual'' DG penalty terms for primary variables, our formulation also has a~special term $\f 1 2 \p_t (\Xi_s-\Xi) \bUx$ in Eq.~\eqref{discrete_general_2} that penalizes the difference between the computed (discontinuous) free surface elevation field and the smoothed (continuous) free surface mesh (Fig.~\ref{Mesh_smoothing}). This term is optional in practical applications but it is indispensable for the proof below to go through. This underscores the importance of consistent treatment of moving free surface geometry. The advantage of including this term is the fact that our stability analysis is not tied to any specific choice of mesh smoothing algorithm.
\end{Remark}
Incorporating our approximations for boundary conditions and the explicit forms of the modified Lax-Friedrichs fluxes as well as summing over all elements, we end up with the following system for the free flow:
\begin{subequations}
\label{discrete_formulation}
\begin{align}
 & \hspace{-9mm}\sum_{\elemx \in I_{e,2D}} \left(\p_t \Xi, \delta \right)_\elemx + A_H(\bUx, \delta) = 0,\label{discrete_formulation_h}\\
 & \hspace{-9mm}\sum_{\elem \in I_e} \left(\p_t \bUx, \bphi \right)_\elem 
+ A_U(\Xi, \bphi)
+ B_U(\bUx, \bU, \bphi)
+ E_U(\Q, \bphi)
+ \!\!\!\sum_{\git \in I_{top}} \!\!\!\lan \f {n_z} 2 \p_t (\Xi_s-\Xi) \bUx, \bphi \ran_\git \!
+ \Lambda_U(\bUx, \bphi) 
= \!\!\sum_{\elem \in I_e}\left(\bF_U, \bphi \right)_\elem, \label{discrete_formulation_u} \\
 & \hspace{-9mm}\sum_{\elem \in I_e} \left(\D^{-1}(\bU)\; \Q, \bpsi \right)_\elem
+ E_Q(\bUx, \bpsi) =  0,\label{discrete_formulation_q}\\
 & \hspace{-9mm}A_H(\bUx, \sigma) 
+ \sum_{\git \in I_{horiz}} \!\!\! \lan \bUd, \jump{\sigma}\ran_\git 
+ \sum_{\git \in I_{top}} \!\! \lan \bU \cdot \bn, \sigma \ran_\git 
- \sum_{\elem \in I_e} \left(W, \p_z \sigma \right)_\elem = 0\label{discrete_formulation_w}
\end{align}
\end{subequations}
with forms $A_H, A_U, B_U, \Lambda_U, E_U, E_Q$ defined as
\begin{align*}
& \hspace{-5mm} A_H(\bUx, \sigma) \coloneqq \sum_{\git \in I_{lat}} \lan \bUxa, \jump{\sigma}\ran_\git
+ \ \sum_{\git \in I_i}  \lan \bUx \cdot \bnx, \sigma \ran_\git\ 
+ \sum_{\git \in I_{bot}} \lan \TU \cdot \bn, \sigma \ran_\git
+ \sum_{\git \in I_o} \lan \buxh \cdot \bnx, \sigma \ran_\git 
-  \sum_{\elem \in I_e} \left(\bUx \cdot \nablax, \sigma \right)_\elem, \nonumber\\
& \hspace{-5mm} A_U(\Xi, \bphi) \coloneqq  \sum_{\git \in I_{lat}} \lan \avg{\Xi}, \jump{\bphi} \ran_\git
+ \ \sum_{\git \in I_{horiz}} \lan \Xi, \jump{\bphi} \ran_\git\ 
+ \ \sum_{\git \in I_i} \lan \hat {\xi} \bnx, \bphi \ran_\git
+ \sum_{\git \in I_o \cup I_{top} \cup I_{bot}} \lan \Xi \bnx, \bphi \ran_\git\ 
- \ \sum_{\elem \in I_e} \left(\Xi \nablax, \bphi \right)_\elem, \nonumber\\
& \hspace{-5mm} B_U(\bUx, \bU, \bphi) \coloneqq 
 \sum_{\git \in I_{lat}} \!\!\! \lan \avg{\bUx \otimes \bU}, \jumpt{\bphi} \ran_\git 
+  \sum_{\git \in I_{horiz}} \!\!\! \lan \bUxa \otimes \bUd, \jumpt{\bphi} \ran_\git
+ \sum_{\git \in I_i} \lan \bUx (\bU \cdot \bng), \bphi \ran_\git\\
 &\hspace{20mm} + \sum_{\git \in I_o} \lan \bUx (\buh \cdot \bng), \bphi \ran_\git
+ \sum_{\git \in I_{top}} \!\!\! \lan \bUx (\bU \cdot \bn), \bphi \ran_\git 
 + \sum_{\git \in I_{bot}} \lan \bUx (\TU \cdot \bn), \bphi \ran_\git
 - \sum_{\elem \in I_e} \left(\bUx(\bU \cdot \nabla), \bphi \right)_\elem, \nonumber\\
& \hspace{-5mm} \Lambda_U(\bUx, \bphi) \coloneqq 
 \sum_{\git \in I_{lat}} \lan \f {\lambda_U} 2 \jumpt{\bUx}, \jumpt{\bphi}\ran_\git 
+  \sum_{\git \in I_i} \lan \f {\lambda_U} 2 \,\left( \bUx - \buxh \right), \bphi \ran_\git,\nonumber\\
& \hspace{-5mm} E_U(\Q, \bphi) \coloneqq  \sum_{\git \in I_{lat} \cup I_{horiz}} \!\!\!\lan \rD \Qa, \jumpt{\bphi} \ran_\git 
+  \sum_{\git \in I_i \cup I_o } \!\!\!\!  \lan \rD \Q \cdot \bn , \bphi \ran_\git 
+  \sum_{\git \in I_{bot}} \!\!\!  \lan C_f \bUx, \bphi \ran_\git
 - \sum_{\elem \in I_e} \!\! \left(\rD \Q \cdot \nabla, \bphi \right)_\elem,\\
&  \hspace{-5mm} E_Q(\bUx, \bpsi) \coloneqq \sum_{\git \in I_{lat} \cup I_{horiz}} \! \lan \bUxa , \jump{\bpsi} \ran_\git 
+  \sum_{\git \in I_i \cup I_o} \! \! \lan \buxh , \bpsi \cdot \bn \ran_\git
+  \sum_{\git \in I_{top} \cup I_{bot}} \lan \bUx, \bpsi \cdot \bn \ran_\git 
- \ \sum_{\elem \in I_e} \left( \bUx, \nabla \cdot \bpsi \right)_\elem. \nonumber
\end{align*}

\subsection{Semi-discrete LDG formulation for the Darcy system}
Analogously to the above section, we formulate the semi-discrete Darcy system

\begin{subequations}
\label{Darcy_compact}
\begin{align}
\label{Darcy_compact_1}
& \sum_{\elem \in \tIe} \left(\p_t \TH, \TestH \right)_\elem 
+ \tilde{E}_{\TH} (\TU, \TestH) 
+ \tilde{\Lambda}_{\TH}(\TH, \TestH) 
= \sum_{\elem \in \tIe} \left(f, \TestH \right)_\elem,\\
\label{Darcy_compact_2}
& \sum_{\elem \in \tIe} \left(\D^{-1}(\TH)\, \TU, \TestU \right)_\elem 
+ \tilde{E}_{\TU} (\TH, \TestU) 
= 0
\end{align}
\end{subequations}
with forms $\tilde{E}_{\TH}, \tilde{E}_{\TU}, \tilde{\Lambda}_{\TH}$ defined as follows:
\begin{align*}
& \hspace{-7mm} \tilde{E}_{\TH} (\TU, \TestH)  \coloneqq  
\sum_{\git \in \tIi} \left\langle \avg{\TU}, \jump{\TestH} \right\rangle_{\git}
+ \sum_{\git \in \tID \cup \tIT } \left\langle \TU \cdot \tilde \bn, \TestH\right\rangle_{\git} 
+ \sum_{\git \in \tIN} \left\langle \hat u_{\tilde n}, \TestH \right\rangle_{\git} 
- \sum_{\elem \in \tIe} \left(\TU \cdot \nabla, \TestH \right)_\elem, \\ 
& \hspace{-7mm} \tilde{E}_{\TU} (\TH, \TestU) \coloneqq 
\sum_{\git \in \tIi} \left\langle \avg{\TH}, \jump{\TestU}\right\rangle_{\git} 
+ \sum_{\git \in \tilde I_{top}} \left\langle \Xi+\f 12 \left(\bUx \cdot \bUx\right), \TestU \cdot \tilde \bn \right\rangle_{\git} 
+ \sum_{\git \in \tIN} \left\langle \TH, \TestU \cdot \tilde \bn \right\rangle_{\git} 
+ \sum_{\git \in \tilde I_D}  \left\langle \hat h, \TestU \cdot \tilde \bn \right\rangle_{\git} 
- \sum_{\elem \in \tIe} \left(\TH, \nabla \cdot \TestU \right)_\elem, \\
& \hspace{-7mm} \tilde{\Lambda}_{\TH}(\TH, \TestH) \coloneqq 
\sum_{\git \in \tIi} \frac{\eta}{\Delta x_{\git}} \left\langle \jump{\TH}, \jump{\TestH} \right\rangle_{\git} 
+ \sum_{\git \in \tID} \frac{\eta}{\Delta x_{\git}} \left\langle \TH - \hat h, \TestH\right\rangle_{\git}.
\end{align*}
The initial state $\TH(0)$ is created by the element-wise $L^2$-projection of $\tH_0$. Here, $\Xi = \Xi(x,y)$ denotes the free surface elevation from the free flow problem.
\section{Discrete energy stability estimate for the coupled system}\label{sec:analysis}
%

%
%

\begin{Theorem}[Discrete stability\label{discrete_stability}]
Let the free surface elevation of the smoothed mesh satisfy $\Xi_s\big|_{\Pi(\p \Omega_{i})}=\hat \xi$, and let $\Xi, \delta \in \IP_{2k}^2(\Pi \elem)$, $\bUx, \bphi \in \IP_k^3(\elem)^2$, $W, \sigma \in \IP_{2k}^3(\elem)$, $\Q, \bpsi \in \IP_k^3(\elem)^{2 \times 3}$, $\TH, \TestH \in \IP_\hat k(\tilde \elem)$, and $\TU, \TestU \in \IP_\bar k^3(\tilde \elem)$ for some $k, \bar k, \hat k \ge 0$, a.e. $t\in[0,T]$, and all $\elem \in \Th(\Omega(t)), \tilde \elem \in \Th(\tO)$. Then scheme~\eqref{discrete_formulation}--\eqref{Darcy_compact} is stable in the following sense:
\begin{align*}
  & \hspace{-8mm} \p_t \left\{ \left\| \Xi \right\|_\Ox^2 
+ \left\| \bUx \right\|^2_{\Omega(T)} 
+ \left\| \TH \right\|^2_\tO \right\}
+ \left\| \sqrt{\D^{-1}(\bU)}\, \Q \right\|^2_{\Ot} 
+ \sum_{\git \in I_{lat}} \left\|\jumpt{\bUx} \right\|^2_\git
+ \bigg\| \sqrt{\tD^{-1}(\TH)}\, \TU \bigg\|^2_{\tO } 
+ \sum_{\git \in \tIi} \frac{\eta}{\Delta x_{\git}} \left\| \jump{\TH} \right\|^2_{\git} 
+ \sum_{\git \in \tID}  \frac{\eta}{\Delta x_{\git}} \left\| \TH \right\|^2_{\git}\\
  & \hspace{-8mm} \qquad \le C(C_t, \Ot, \tO, \D, \tD, \eta, \bF_U, \tilde f, z_b, \hat{\xi}, \buxh, \hat u_{\tilde n}, \Delta x).
 \end{align*}
\end{Theorem}
\noindent
{\em Proof}:
We start with the stability estimate for Darcy flow. Choosing $\TestH = \TH$, $\TestU = \TU$ and adding \eqref{Darcy_compact_1}, \eqref{Darcy_compact_2} gives
\begin{equation} \label{Darcy_stability_1}
\sum_{\elem \in \tIe} \left(\p_t \TH, \TH \right)_\elem 
+ \tilde{E}_{\TH} (\TU, \TH) 
+ \tilde{\Lambda}_{\TH}(\TH, \TH) 
+ \sum_{\elem \in \tIe} \left(\tD^{-1}(\TH)\, \TU, \TU \right)_\elem 
+ \tilde{E}_{\TU} (\TH, \TU) 
= \sum_{\elem \in \tIe} \left(f, \TH \right)_\elem.
\end{equation}
Integration by parts of the element integral term in $\tilde{E}_{\TU}$ and the use of \eqref{jump1} leads to
\[
\tilde{E}_{\TH} (\TU, \TH) + \tilde{E}_{\TU} (\TH, \TU) 
= \sum_{\git \in \tIT} \left\langle \Xi +\f 12 (\bUx \cdot \bUx), \TU \cdot \tilde \bn \right\rangle_{\git}
+ \sum_{\git \in \tID} \left\langle \hat h, \TU \cdot \tilde \bn \right\rangle_{\git} 
+ \sum_{\git \in \tIN} \left\langle \hat u_{\tilde n}, \TH \right\rangle_{\git}.
\]
Substituting the above expression into~\eqref{Darcy_stability_1} and splitting the penalty terms results in
\begin{align*}
 & \hspace{-8mm} \frac{1}{2} \p_t \left\| \TH \right\|^2_{\tO }
+ \left\| \sqrt{\tD^{-1}(\TH)}\, \TU \right\|^2_{\tO}
+ \sum_{\git \in \tIi} \frac{\eta}{\Delta x_{\git}} \left\| \jump{\TH} \right\|^2_{\git} 
+ \sum_{\git \in \tID} \frac{\eta}{\Delta x_{\git}} \left\| \TH \right\|^2_{\git}
+ \sum_{\git \in \tIT} \left\langle \Xi +\f 12 (\bUx \cdot \bUx), \TU \cdot \tilde \bn \right\rangle_{\git} \\
 & \hspace{-8mm}\quad = \un{\left(\tilde f, \TH\right)_{\tO}}_{\tilde \Upsilon_1}
+ \un{\sum_{\git \in \tID} \frac{\eta}{\Delta x_{\git}} \left\langle\hat h, \TH \right\rangle_{\git} }_{\tilde \Upsilon_2}
- \un{\sum_{\git \in \tID} \left\langle \hat h, \TU \cdot \tilde \bn \right\rangle_{\git} }_{\tilde \Upsilon_3}
- \un{\sum_{\git \in \tIN} \left\langle \hat u_{\tilde n}, \TH \right\rangle_{\git}}_{\tilde \Upsilon_4}
\end{align*}
Now we estimate terms $\tilde \Upsilon_1$--$\tilde \Upsilon_4$ using Young's inequality, uniform bounds on $\tD(\TH)$, and the auxiliary results from Sec.~\ref{sec:tools} (also see \cite[p. 1382 - 1383]{RuppKnabner2017} presenting similar estimates in greater detail).
\begin{align*}
\hspace{-8mm}|\tilde{\Upsilon}_1| & 
\;\le\; \f 12 \|\tilde f\|^2_{\tO} + \f 12 \|\TH\|^2_{\tO},\\
\hspace{-8mm}|\tilde{\Upsilon}_2| & 
\;\le\; \sum_{\git \in \tID} \frac{\eta}{2 \Delta x_{\git}} \| \hat h \|^2_{\git}
+ \sum_{\git \in \tID} \frac{\eta}{2 \Delta x_{\git}} \| \TH \|^2_{\git},\\
\hspace{-8mm}|\tilde{\Upsilon}_3| & \;\le\;
\sum_{\git \in \tID} \frac{C_t C_{D}}{2 \Delta x_{\git}} \| \hat h \|^2_{\git}
+ \sum_{\git \in \tID} \frac{\Delta x_{\git}}{2 C_t C_{D}} \| \TU \|^2_{\git}
\;\le\; C(C_t, \tD) \sum_{\git \in \tID} \Delta x_{\git}^{-1} \| \hat h \|^2_{\git}
+ \frac 12  \left\| \sqrt{\tD^{-1}(\TH)}\, \TU \right\|^2_{\tO},\\
\hspace{-8mm}|\tilde{\Upsilon}_4| & 
\;\le\; \sum_{\git \in \tIN} \frac{C_t}{2 \Delta x_{\git}} \| \hat u_{\tilde n} \|^2_{\git}
+ \sum_{\git \in \tIN} \frac{\Delta x_{\git}}{2 C_t} \| \TH \|^2_{\git}
\;\le\; C(C_t) \sum_{\git \in \tIN} \Delta x_{\git}^{-1} \| \hat u_{\tilde n} \|^2_{\git}
+ \f 12 \|\TH\|^2_{\tO}.
\end{align*}
Using the above estimates and noting that $\tilde \bn = - \bn$ on $\p \tO_{top}$, we obtain
\begin{align}\label{EQ:Darcy:estimate}
 & \f 12 \p_t \left\| \TH \right\|^2_{\tO} 
+ \f 12 \left\| \sqrt{\tD^{-1}(\TH)}\, \TU \right\|^2_{\tO} 
+ \sum_{\git \in \tIi} \frac{\eta}{\Delta x_{\git}} \left\| \jump{\TH} \right\|^2_{\git} 
+ \sum_{\git \in \tID} \frac{\eta}{2 \Delta x_{\git}} \left\| \TH \right\|^2_{\git}
- \sum_{\git \in \tIT} \left\langle \Xi+\f 12 (\bUx \cdot \bUx), \TU \cdot \bn \right\rangle_{\git} \notag\\
 & \quad \le \| \TH \|^2_\tO 
+ \f 12 \| \tilde f \|^2_{\tO} 
+ C(C_t, \tD, \eta) \sum_{\git \in \tID} \Delta x_{\git}^{-1} \| \hat h \|^2_{\git}
+ C(C_t) \sum_{\git \in \tIN} \Delta x_{\git}^{-1} \| \hat u_{\tilde n} \|^2_{\git}.
\end{align}
\noindent
Turning to the analysis of the free surface flow sub-system, we set $\delta=\Xi, \bphi=\bUx, \bpsi =\Q$ and add Eqs.~\rf{discrete_formulation_h}--\rf{discrete_formulation_q}
\begin{align}
\label{discrete_stability_1}
& \sum_{\elemx \in I_{e,2D}} \left(\p_t \Xi, \Xi \right)_\elemx 
 +  A_H(\bUx, \Xi) 
 + \sum_{\elem \in I_e} \left(\p_t \bUx, \bUx \right)_\elem 
+ A_U(\Xi, \bUx) 
+ B_U(\bUx, \bU, \bUx) 
+ E_U(\Q, \bUx) 
+ \Lambda_U(\bUx, \bUx)  
\nonumber\\
& \quad + \sum_{\git \in I_{top}} \!\!\! \lan \f {n_z} 2 \p_t (\Xi_s-\Xi) \bUx, \bUx \ran_\git
 + \sum_{\elem \in I_e} \!\!\left(\D^{-1}(\bU)\;\Q, \Q \right)_\elem 
+   E_Q(\bUx, \Q)  
=  \sum_{\elem \in I_e}\!\!\left(\bF_U, \bUx \right)_\elem. 
\end{align}
First, we deal with  terms containing $\Xi$. Integration by parts and \rf{jump1} produce
\[
A_H(\bUx, \Xi) +  A_U(\Xi, \bUx)  
= \sum_{\git \in I_{bot}}  \lan \TU \cdot \bn, \Xi \ran_\git
+ \sum_{\git \in I_i} \lan \hat{\xi} \bnx, \bUx \ran_\git
+ \sum_{\git \in I_o}  \lan \buxh \cdot \bnx, \Xi \ran_\git.
\]
The step dealing with the non-linear advective terms in the momentum equation is the crucial and, 
at the same time, the most technically involved step of the proof and thus will be presented in 
greater detail.
First, note that 
\[
- \sum_{\elem \in I_e} \left(\bUx(\bU \cdot \nabla), \bUx \right)_\elem
\ = \ - \sum_{\elem \in I_e} \left(\bU, \f 1 2 \nabla \left(\bUx \cdot \bUx \right) \right)_\elem.
\]
Setting $\sigma = \f 1 2 \left(\bUx \cdot \bUx \right)$ in~\rf{discrete_formulation_w} (recalling that its test space contains products of elements from the test space of \rf{discrete_formulation_u}!), we replace the above term in the definition of $B_U$ with the boundary integral terms resulting from~\rf{discrete_formulation_w}. 
\begin{align*}
\hspace{-8mm}B_U(\bUx, \bU, \bUx) 
&=  \sum_{\git \in I_{lat}} \bigg\{ \un{ \lan \avg{\bUx \otimes \bU}, \jumpt{\bUx} \ran_\git }_{\Theta_1}
- \f 1 2 \lan \bUxa, \jump{ \bUx \cdot \bUx }\ran_\git  \bigg\}
+ \!\!\sum_{\git \in I_{horiz}} \!\! \bigg\{\un{\lan \bUxa \otimes \bUd, \jumpt{\bUx} \ran_\git}_{\Theta_2}
- \f 1 2 \lan \bUd, \jump{\bUx \cdot \bUx } \ran_\git \bigg\}\\
\hspace{-8mm}& \quad + \f 1 2 \sum_{\git \in I_o} \lan \bUx\, (\buxh \cdot \bnx), \bUx \ran_\git
 + \f 1 2 \sum_{\git \in I_{top} \cup I_i} \lan \bUx\, (\bU \cdot \bn), \bUx \ran_\git 
+  \f 1 2 \sum_{\git \in I_{bot}} \lan \bUx\, (\TU \cdot \bn), \bUx \ran_\git.
\end{align*}
Using \rf{jump1}, \rf{jump2} and noting $(\bU \cdot \bn)|_{\git} = (\bUx \cdot \bnx)|_{\git}, \forall \git \in I_{lat}$ we find
\begin{align*}
 \Theta_1 &= \lan \avg{\bUx} \otimes \avg{\bU} + \f 1 4 \jumpt{\bUx} \jumpt{\bU}^T, \jumpt{\bUx} \ran_\git
= \lan \avg{\bUx} \otimes \avg{\bU} + \f 1 4 \left( \bUx^+ - \bUx^-\right) \otimes \left( \bU^+ - \bU^-\right), \bUx^+ \otimes \bn^+ + \bUx^- \otimes \bn^- \ran_\git\\
&= \lan \avg{\bUx} \left( \avg{\bU} \cdot \bn^+ \right) + \f 1 4\left( \bUx^+ - \bUx^- \right) \jump{\bU}, \bUx^+ - \bUx^- \ran_\git
= \f 1 2 \lan \avg{\bUx}, \jump{\bUx \cdot \bUx} \ran_\git
+ \f 1 4 \lan \jumpt{\bUx} \jump{\bU}, \jumpt{\bUx} \ran_\git.
\end{align*}
In a~similar manner, we obtain $\Theta_2 = \f 1 2 \lan \bUd, \jump{\bUx \cdot \bUx } \ran_\git$; this gives us
\begin{equation}\label{B_U}
B_U(\bUx, \bU, \bUx) 
= \f 1 4 \sum_{\git \in I_{lat}} \lan \jumpt{\bUx} \jump{\bU}, \jumpt{\bUx} \ran_\git
+ \f 1 2 \bigg\{ \sum_{\git \in I_o} \lan \bUx\, (\buxh \cdot \bnx), \bUx \ran_\git
+ \sum_{\git \in I_{top} \cup I_i} \lan \bUx\, (\bU \cdot \bn), \bUx \ran_\git 
+  \sum_{\git \in I_{bot}} \lan \bUx\, (\TU \cdot \bn), \bUx \ran_\git \bigg\}.
\end{equation}

The movement of the free surface is accounted for via the mesh penalty term. 
Here we use Eqs. \rf{discrete_formulation_h} and \rf{discrete_formulation_w}, 
once again taking advantage of the higher order test spaces in them. 
Noting that for any free surface boundary face $\gamma$ in our smoothed mesh 
$\int_{\gamma} n_z f(x, y, z) ds = \int_{\Pi(\gamma)} f(x, y, \Xi_s) dx dy$
and applying Leibniz' Rule we proceed as follows:
\begin{align*}
&\hspace{-5mm} \sum_{\elem \in I_e} \left(\p_t \bUx, \bUx \right)_\elem
 + \ \sum_{\git \in I_{top}} \lan \f {n_z} 2 \p_t (\Xi_s-\Xi) \bUx, \bUx \ran_\git 
= \  \f 1 2 \sum_{\elemx \in I_{e,2D}} \!\!\! \left( \int_{z_{bot}}^{\Xi_s(t)} 
\p_t |\bUx|^2 dz, 1 \right)_\elemx
+ \f 1 2 \sum_{\elemx \in I_{e,2D}} \!\!\! \left( \p_t \Xi_s - 
\p_t \Xi, |\bUx(\Xi_s)|^2 \right)_\elemx  \nonumber\\
& \hspace{-5mm} \quad = \  \f 1 2 \sum_{\elem \in I_e} \p_t \left\| \bUx \right\|^2_\elem 
\ - \ \f 1 2 \sum_{\elemx \in I_{e,2D}} \left(\p_t \Xi, |\bUx(\Xi_s)|^2 \right)_\elemx
=\  \f 1 2 \sum_{\elem \in I_e} \p_t \left\| \bUx \right\|^2_\elem 
\ - \ \f 1 2 \sum_{\git \in I_{top}} \lan \bU \cdot \bn, \bUx \cdot \bUx \ran_\git.
\end{align*}
The last equality follows by setting $\delta = \sigma = |\bUx(\Xi_s)|^2$ in
\rf{discrete_formulation_h} and
\rf{discrete_formulation_w}, respectively, and subtracting the latter from the former.
The last term in the expression above cancels a corresponding term in the estimate~\eqref{B_U} for $B_U$.

For the diffusion terms, the divergence theorem and \rf{jump1} give us
\[
E_U(\Q, \bUx) + E_Q(\bUx, \Q)
= \sum_{\git \in I_i \cup I_o} \lan \buxh \rD, \Q \cdot \bn \ran_\git
+ \sum_{\git \in I_{bot}} \lan C_f \bUx , \bUx \ran_\git. 
\]
Substituting the results of the above simplifications into \rf{discrete_stability_1} we obtain
\begin{align}
\label{discrete_stability_2}
&\f 1 2 \p_t \left\| \Xi \right\|_\Ox^2
+ \f 1 2 \p_t \left\| \bUx \right\|^2_{\Omega(t)} 
+ \left\| \sqrt{\D^{-1}(\bU)}\, \Q \right\|^2_{\Omega(t)} 
+ \sum_{\git \in I_{lat}} \lan \f {\lambda_U} 2 \jumpt{\bUx}, \jumpt{\bUx} \ran_\git 
+ \sum_{\git \in I_i} \lan \f {\lambda_U} 2 \bUx, \bUx \ran_\git \nonumber\\
& \qquad + \f 1 2  \sum_{\git \in I_{bot}} \lan \bUx  (\TU \cdot \bn), \bUx \ran_\git
+  \un{\sum_{\git \in I_{bot}}  \lan C_f \bUx, \bUx \ran_\git }_{\ge 0} 
+ \un{\f 1 2 \sum_{\git \in I_o} \lan \bUx (\buxh \cdot \bnx), \bUx \ran_\git }_{\ge 0} \nonumber\\
& \quad =  -\un{\f 1 4  \sum_{\git \in I_{lat}} \lan \jumpt{\bUx}\, \jump{\bU}, \jumpt{\bUx} \ran_\git }_{\Upsilon_1}
+ \un{\sum_{\elem \in I_e} \left(\bF_U(t), \bUx \right)_\elem }_{\Upsilon_2} 
 -  \un{ \sum_{\git \in I_i} \lan \hat{\xi} \bnx, \bUx \ran_\git }_{\Upsilon_3}  
- \un{\sum_{\git \in I_o} \lan \buxh \cdot \bnx, \Xi \ran_\git }_{\Upsilon_4} 
+  \un{\sum_{\git \in I_i} \lan \f {\lambda_U} 2 \buxh, \bUx \ran_\git }_{\Upsilon_5} \nonumber\\
& \qquad - \un{\sum_{\git \in I_i \cup I_o} \lan \buxh \rD, \Q \cdot \bn \ran_\git }_{\Upsilon_6}
 - \un{\f 1 2 \sum_{\git \in I_i} \lan \bUx (\bU \cdot \bn), \bUx \ran_\git }_{\Upsilon_7}
- \sum_{\git \in I_{bot}} \lan \TU \cdot \bn, \Xi \ran_\git.
\end{align}
In the remainder of the proof, we estimate terms $\Upsilon_1$--$\Upsilon_7$ relying on 
Young's and Cauchy-Schwarz' inequalities, properties of $\Xi_s$ and $\D$, and results from Sec.~\ref{sec:tools}.
\begin{align*}
\hspace{-8mm}|\Upsilon_1| & \le \f 1 2 \sum_{\git \in I_{lat}} \lan \jumpt{\bUx} \avg{|\bU \cdot \bng|}, \jumpt{\bUx} \ran_\git, \qquad \text{(cf. the definition of $\lambda_U$)}\\
\hspace{-8mm}|\Upsilon_2| & \le \f 1 4 \left\| \bF_U(t) \right\|_{\Omega(t)}^2
 + \; \left\| \bUx \right\|_{\Omega(t)}^2,\\
\hspace{-8mm}|\Upsilon_3| & \le \f 1 {4 \alpha} \sum_{\git \in I_i} \| \hat{\xi} \|^2_\git
 + \alpha \sum_{\git \in I_i} \| \bUx \|^2_\git,
\; \mbox{where} \;  \f 1 {4 \alpha} \sum_{\git \in I_i} \| \hat{\xi} \|^2_\git
= \f 1 {4 \alpha} \lan \Xi_s-z_b, \hat{\xi}^2 \ran_{\Pi(\p \Omega_i)} 
\le \f 1 {4 \alpha} \left( \| \hat{\xi} \|^3_{\Pi(\p \Omega_i)}
+ \|z_b\|_{L^\infty(\p \Omega_i)}\| \hat{\xi} \|^2_{\Pi(\p \Omega_i)} \right), \\
\hspace{-8mm}|\Upsilon_4| & = \lan \Xi, \int^{\Xi_s}_{z_b} \buxh \cdot \bnx \, dz \ran_{\Pi(\p \Omega_o)}
\le \f{\Delta x}{C_t} \| \Xi\|^2_{\Pi(\p \Omega_o)}
+ \f{C_t}{4 \Delta x} \lan 1, \left(\int^{\Xi_s}_{z_b} \buxh \cdot \bnx \, dz \right)^2 \ran_{\Pi(\p \Omega_o)}
\le \| \Xi\|^2_{\Ox}
+ C(C_t, \Delta x, \p \Omega_o)\, \| \buxh \|^2_{\p \Omega_o},\\
\hspace{-8mm}|\Upsilon_5| & \le \f 1 {8 \beta} \;\sum_{\git \in I_i} \lan \lambda_U \buxh, \buxh \ran_\git
+\; \beta \sum_{\git \in I_i} \lan \lambda_U \bUx, \bUx \ran_\git
\; \le \; \sum_{\git \in I_i} \lan \lambda_U, \lambda_U \ran_\git
+ \;\f 1 {32 \beta} \sum_{\git \in I_i}  \lan |\buxh|^2, |\buxh|^2 \ran_\git
 + \; \beta \sum_{\git \in I_i} \lan \lambda_U \bUx, \bUx \ran_\git \nonumber\\
\hspace{-8mm}& \le  \alpha \sum_{\git \in I_i} \| \bUx \|^2_\git
+ \f 1 {4 \alpha} \sum_{\git \in I_i} \lan 1, 1 \ran_\git
+ \f 1 {32 \beta} \sum_{\git \in I_i} \|\buxh\|^4_\git
 +  \beta \sum_{\git \in I_i} \lan \lambda_U \bUx, \bUx \ran_\git 
\; \mbox{with}\quad \f 1 {4 \alpha} \sum_{\git \in I_i} \lan 1, 1 \ran_\git 
= \f 1 {4 \alpha} \lan 1, \hat{\xi} -z_b \ran_{\Pi(\p \Omega_i)}\\
\hspace{-8mm}|\Upsilon_6| & \le \f {C_t^2 C_{D}} {2} \sum_{\git \in I_i \cup I_o} \Delta x_{\git}^{-1} \| \buxh\|^2_\git
\;+ \; \f 1 {2 C_t^2 C_{D}} \sum_{\git \in I_i \cup I_o} \Delta x_{\git} \|\Q \|^2_\git 
 \;\le\; C\left(C_t, \D\right) \sum_{\git \in I_i \cup I_o} \Delta x_{\git}^{-1} \| \buxh\|^2_\git
+\; \f 1 2 \left\|\sqrt{\D^{-1}(\bU)}\, \Q \right\|^2_{\Omega(t)},
\end{align*}
where $0 < \alpha, \beta < 1$ are some parameters that will be determined later.

Collecting the terms containing $\bUx$ on the inflow faces, namely $\Upsilon_7$ and the corresponding 
terms in the estimates for $\Upsilon_3$ and $\Upsilon_5$, we use the penalty term to estimate on the left hand side of
\eqref{discrete_stability_2}
\[
\sum_{\git \in I_i} \lan \f {\lambda_U} 2 \bUx, \bUx \ran_\git
 - \; 2 \alpha \sum_{\git \in I_i} \| \bUx \|^2_\git
- \; \beta \sum_{\git \in I_i}  \lan \lambda_U \bUx, \bUx \ran_\git  
 - \; \f 1 2 \sum_{\git \in I_i}  \lan \bUx \left|\bU \cdot \bn \right|, \bUx \ran_\git 
 \ge 0,
\]
which by Lemma~\ref{properties_lambda} can be shown to hold for the following choices of $\alpha$ and $\beta$:
\[
0 \; < \; \beta \; \le \; \f 1 {2 \sqrt{2} + 2}, \qquad 0 \; < \; \alpha \;\le\; \f {1/2 - \beta} {2\sqrt{2}}. 
\]
Substituting the estimates above into~\rf{discrete_stability_2} we obtain the following inequality:
\begin{align}\label{discrete_stability_3}
 & \f 1 2 \p_t \left\| \Xi \right\|_\Ox^2 
+ \f 1 2 \p_t  \left\| \bUx \right\|^2_{\Omega(t)} 
+ \f 1 2 \left\| \sqrt{\D^{-1}(\bU)}\, \Q(t) \right\|^2_{\Omega(t)} 
+ \f 1 2 \sum_{\git \in I_{lat}}\!\!\! \lan \jumpt{\bUx}, \jumpt{\bUx} \ran_\git
+ \sum_{\git \in I_{bot}} \lan \Xi + \f 1 2 ( \bUx \cdot \bUx ), \TU \cdot \bn  \ran_\git 
\nonumber\\
 & \qquad \le \| \Xi\|^2_{\Ox}
+ \left\| \bUx \right\|_{\Omega(t)}^2 
+ \; C(C_t, \Ot, \bF_U, z_b, \hat{\xi}, \buxh, \Delta x).
\end{align}

The claim of our theorem follows by adding~\eqref{EQ:Darcy:estimate} to~\eqref{discrete_stability_3}.

\framebox[\width]{\Huge{ }}\\

\section{Numerical results}\label{sec:numerical}
The numerical implementation is based on our FESTUNG framework \cite{FrankRAK2015,ReuterAWFK2016,JaustRASK2018} and, specifically, utilizes the setup detailed in the companion paper~\cite{ReuterRAK2018}.
We choose a two-dimensional (in a~vertical $xz$-slice) computational domain~$\Omega(t)\cup\tO \subset\IR^2$ with~$\Omega(t) \coloneqq (0,100)\times(z_b,\xi(t))$, $\tO \coloneqq (0,100)\times(-5,z_b)$, time interval~$J=(0,10)$, and a~sloped interface between free flow and subsurface domains~$z_b(x^1) \coloneqq 0.005 x^1$, which has a~constant normal vector~$\vec{\nu} = \pm 1/\sqrt{1+0.005^2} \,[-0.005, 1]^T$.
For a~given free surface elevation~$\xi$ and horizontal velocity~$u$, one can derive matching analytical functions for~$\tH$ using interface condition~\eqref{transition_bc_2} and for~$w$ using continuity equation~\eqref{cont_eq} and interface condition~\eqref{transition_bc_1}.
Instead of \eqref{hydrostatic_surface_bc}, we use here non-homogeneous boundary conditions at the free surface to have more freedom in our choice for~$u$ resulting in the following analytical solution
\begin{align*}
\xi(t,x) &\coloneqq\; 
  5 + 0.003\,\sin(0.08\, x + 0.08\, t) \,,\\
u(t,\vec{x}) &\coloneqq\; 
  r(t,x) \big( \cos(0.1\, z) - \cos\left(0.1\, z_b(x)\right) \big)\,,\\
w(t,\vec{x}) &\coloneqq\; 
  n(t,\vec{x}) + \varepsilon(t,x) \,,\\
\tH(t,\vec{x}) &\coloneqq\;
  \xi(t,x) + \big(\sin(0.3\, z) - \sin\left(0.3\, z_b(x)\right)\big)\, m(t,x)
\end{align*}
with diffusion coefficients~$\D \coloneqq 0.05 \,I$, $\tD \coloneqq 0.01\, I$. $n(t,\vec{x})$ is chosen so that~$\partial_{x} u + \partial_{z} w = 0$ in~$\Omega(t)$:
\begin{equation*}
n(t,\vec{x}) \,\coloneqq\; 
  -\partial_{x} r(t,x) \left( \frac{1}{0.1} \sin(0.1\, z) - z \cdot \cos\left(0.1 \,z_b(x) \right)  \right)  
 - 0.1 \cdot 0.005\, \cdot r(t,x) \, z \cdot \sin\left(0.1\, z_b(x)\right)\,,
\end{equation*}
and $\varepsilon(t,x)$ shifts~$w$ to fulfill coupling condition~\eqref{transition_bc_2}, i.e., 
\begin{equation*}
\varepsilon(t,x) \,\coloneqq\; 0.01 \left( 0.005 \, \partial_{x}\tH\left(t,x,z_b(x)\right) - \partial_{z} \tH\left(t,x,z_b(x)\right) \right) - n\left(t,x,z_b(x)\right)\,.
\end{equation*}
Functions~$r(t,x),m(t,x)$ are used to increase the spatial variability in $x$-direction and to introduce a~time dependency.
Here, we use
\begin{equation*}
r(t,x)\,\coloneqq\; \sin(0.07\, x + 0.4\, t) 
\quad\text{ and }\quad
m(t,x)\,\coloneqq\; \cos(0.07\, x + 0.07\, t) \,.
\end{equation*}
We prescribe Dirichlet boundary conditions for all unknowns and derive boundary data, right hand side functions, and initial data from the analytical solution.
Using this setup, we compute the solution for a~sequence of increasingly finer meshes with element sizes~$\Delta x_j$ and evaluate errors and estimated orders of convergence for any function~$c_{\Delta}$ by
\begin{equation*}
\mathrm{Err}(c) \,\coloneqq\; \|c_{\Delta_{j-1}} - c\|_{L^2(\Omega)}\,, \qquad\qquad
\mathrm{EOC}(c) \,\coloneqq\; \ln \left(\frac{\|c_{\Delta_{j-1}} - c\|_{L^2(\Omega)}}{\|c_{\Delta_{j}} - c\|_{L^2(\Omega)}} \right)\Bigg/ \ln \left(\frac{\Delta x_{j-1}}{\Delta x_j}\right)
\end{equation*}
and list those in~Table~\ref{tab:conv:coupled}. Following our analysis, we use approximations of polynomial order $2p$ for $h$ and $w$, whereas all other unknowns are approximated with order $p$.

\begin{table}[!ht]
\small
\setlength{\tabcolsep}{4pt}
\renewcommand{\arraystretch}{1.1}
\begin{tabular}{cccccccccccccc}
\toprule
$p$ & $j$ & Err($\xi$) & $\EOC{\xi}$ & Err($u$) & $\EOC{u}$ & Err($w$) & $\EOC{w}$ & Err($\tH$) & $\EOC{\tH}$ & Err($\tilde u$) & $\EOC{\tilde u}$ & Err($\tilde w$) & $\EOC{\tilde w}$\\
\bottomrule
  & 0 & 2.47e-01 & ---  & 9.63e-01 & ---  & 2.40e-01 & ---  & 4.60e+00 & ---  & 3.95e-01 & ---  & 1.47e+00 & ---  \\
  & 1 & 5.52e-02 & 2.16 & 2.16e-01 & 2.16 & 1.17e-01 & 1.03 & 1.53e+00 & 1.59 & 2.94e-01 & 0.43 & 7.65e-01 & 0.94 \\
1 & 2 & 1.43e-02 & 1.95 & 5.62e-02 & 1.94 & 5.85e-02 & 1.00 & 4.08e-01 & 1.90 & 2.12e-01 & 0.47 & 3.96e-01 & 0.95 \\
  & 3 & 3.59e-03 & 1.99 & 1.62e-02 & 1.80 & 2.85e-02 & 1.04 & 9.83e-02 & 2.05 & 1.09e-01 & 0.95 & 1.88e-01 & 1.08 \\
  & 4 & 9.02e-04 & 1.99 & 5.90e-03 & 1.46 & 1.41e-02 & 1.01 & 2.33e-02 & 2.08 & 5.42e-02 & 1.01 & 9.27e-02 & 1.02 \\
\midrule
  & 0 & 1.38e-01 & ---  & 1.25e-01 & ---  & 4.35e-02 & ---  & 1.60e+00 & ---  & 2.82e-01 & ---  & 5.15e-01 & ---  \\
  & 1 & 4.63e-02 & 1.57 & 3.49e-02 & 1.84 & 1.96e-02 & 1.15 & 2.24e-01 & 2.84 & 7.43e-02 & 1.93 & 1.64e-01 & 1.65 \\
2 & 2 & 9.02e-03 & 2.36 & 4.98e-03 & 2.81 & 4.44e-03 & 2.14 & 3.89e-02 & 2.52 & 2.28e-02 & 1.70 & 4.39e-02 & 1.90 \\
  & 3 & 2.01e-03 & 2.17 & 7.02e-04 & 2.83 & 1.51e-03 & 1.56 & 5.41e-03 & 2.85 & 5.83e-03 & 1.97 & 8.96e-03 & 2.29 \\
  & 4 & 4.69e-04 & 2.10 & 1.32e-04 & 2.41 & 6.81e-04 & 1.15 & 7.04e-04 & 2.94 & 1.47e-03 & 1.99 & 1.81e-03 & 2.31 \\
\bottomrule
\end{tabular}
\caption{$L^2(\Omega)$-errors and estimated orders of convergence (EOC) for the coupled problem. On the $j$th refinement level, we used $2^{j+1} \times 2^j$ elements and time step~$\Delta \tilde{t} = \frac{1}{5} \cdot 2^{-p} \cdot 4^{-j}$ for the subsurface problem and~$\Delta t = \frac{1}{50} \cdot 2^{-p} \cdot 4^{-j}$ for the free flow problem.}
\label{tab:conv:coupled}
\end{table}

\section{Conclusions}
Our stability analysis for the discrete formulation of the coupled hydrostatic/Darcy system motivated our choice of the transition condition for the hydrostatic pressure/hydraulic head. This transition condition includes a~special form of dynamic pressure -- modified to suit the specifics of the hydrostatic model used in the free surface flow system. Further investigations (involving numerical studies and possibly also experimental validations) of this interface condition might be needed to substantiate the physical validity of our choice.

\bibliography{bibliography}   
\bibliographystyle{elsarticle-num}

\end{document}

%% file: Commands.tex
\renewcommand{\bar}[1]{{\overline{#1}}}
\renewcommand{\hat}[1]{{\widehat{#1}}}
\renewcommand{\tilde}[1]{{\widetilde{#1}}}

\newcommand*{\llbrace}{\left\lbrace\!\left\vert}\newcommand*{\rrbrace}{\right\vert\!\right\rbrace}
\newcommand*{\avg}[1]{\llbrace{#1}\rrbrace}                    
\newcommand*{\jump}[1]{\llbracket{#1}\rrbracket} 
\newcommand*{\jumpt}[1]{\,\llbracket\!\!\!\llbracket{#1}\rrbracket\!\!\!\rrbracket\,} 

\newcommand{\IR}{\mathbb R}

\newcommand{\IP}{\mathbb P}

\newcommand{\F}{\mathcal F}

\renewcommand{\vec}{\mathbf}
\newcommand{\gvec}{\boldsymbol}

\newcommand {\ba} {\begin{array}}
\newcommand {\ea} {\end{array}}
\newcommand {\be} {\begin{equation}}
\newcommand {\ee} {\end{equation}}
\newcommand {\bea} {\begin{eqnarray}}
\newcommand {\eea} {\end{eqnarray}}
\newcommand {\beas} {\begin{eqnarray*}}
\newcommand {\eeas} {\end{eqnarray*}}
\newcommand {\lan} {\left\langle}
\newcommand {\ran} {\right\rangle}
\newcommand {\f} {\frac}
\newcommand {\p} {\partial}
\newcommand {\un} {\underbrace}
\newcommand{\bn} {{\bf n}}
\newcommand{\bng} {{\bf n}}
\newcommand{\bnx} {{\bf n}_{2}}
\newcommand{\rf}[1]{(\ref{#1})}
\newcommand{\bu}{{\bf u}}
\newcommand{\bU}{{\bf U}}

\newcommand{\bUx}{{\bf U}_{2}}
\newcommand{\bUxa}{\avg{\bf U_{2}}}
\newcommand{\bUd}{{\bf U}^\downarrow}

\newcommand{\bux}{{\bf u}_{2}}
\newcommand{\buh}{\hat{\bf u}}
\newcommand{\buxh}{\hat{\bf u}_{2}}

\newcommand{\nablax}{\nabla_{2}}
\newcommand{\bF} {{\bf F}}
\newcommand{\D}{{\mathcal{D}}}
\newcommand{\rD}{}
\newcommand{\q}{{\scriptstyle{\mathcal{Q}}}}
\newcommand{\Q}{{\mathcal{Q}}}
\newcommand{\Qa}{\avg{\mathcal{Q}}}

\newcommand{\elem}{{\mathcal K}}
\newcommand{\elemx}{{\mathcal K_2}}

\newcommand{\Ot}{{\Omega(t)}}

\newcommand{\Ox}{{\Omega_{2}}}

\newcommand{\git}{{\gamma}}
\newcommand{\bphi}{{\gvec \varphi}}
\newcommand{\bpsi}{{\it \Psi}}

\newcommand{\Th}{{{\mathcal T}_{\Delta}}}
\newcommand{\bR}{{\rm I}\!{\rm R}}      

\newcommand{\tO}{\tilde{\Omega}}
\newcommand{\tU}{\tilde{\vec u}}
\newcommand{\tH}{\tilde h}
\newcommand{\tD}{\tilde \D}
\newcommand{\testU}{\tilde{\bphi}}
\newcommand{\testH}{\tilde{\delta}}
\newcommand{\tIe}{\tilde I_e}
\newcommand{\tIi}{\tilde I_{int}}
\newcommand{\tID}{\tilde I_{D}}
\newcommand{\tIT}{\tilde I_{top}}
\newcommand{\tIN}{\tilde I_N}
\newcommand{\TU}{\tilde{\vec U}}
\renewcommand{\TH}{\tilde H}
\newcommand{\TestU}{\testU}
\newcommand{\TestH}{\testH}

\newcommand{\EOC}[1]{{\ensuremath{\text{\scriptsize EOC}(#1)}}}